\documentclass[review]{elsarticle}

\usepackage{hyperref}
\usepackage{amssymb}
\usepackage{amsmath,rotating,subfig}
\usepackage{multirow}

\journal{Journal of \LaTeX\ Templates}

\begin{document}

\begin{frontmatter}

\title{A Relax-and-Decomposition Algorithm for a p-Robust Hub Location Problem\tnoteref{mytitlenote}}

\author[mymainaddress]{Saeid Abbasi-Parizi}
\ead{saeed\_abasi@aut.ac.ir}

\author[mysecondaryaddress]{Mahdi Bashiri\corref{mycorrespondingauthor}}
\cortext[mycorrespondingauthor]{Corresponding author}
\ead{bashiri@shahed.ac.ir}

\author[thirdaddress]{Andrew Eberhard}
\ead{andy.eb@rmit.edu.au}

\address[mymainaddress]{Department of Industrial Engineering \& Management Systems, Amirkabir University of Technology, Tehran, Iran}
\address[mysecondaryaddress]{Department of Industrial Engineering, Shahed University, Tehran, Iran}
\address[thirdaddress]{Mathematical and Geospatial Sciences, School of Science, RMIT university, Melbourne 3000, Australia}

\begin{abstract}
In this paper, a non-linear $p$-robust hub location problem is extended to a risky environment where augmented chance constraint with a min-max regret form is employed to consider network risk as one of the objectives. The model considers risk factors such as security, air pollution and congestion to design the robust hub network. A Monte-Carlo simulation based algorithm, namely, a sample average approximation scheme is applied to select a set of efficient scenarios. The problem is then solved using a novel relax-and-decomposition heuristic based on the coupling of an accelerated Benders decomposition with a Lagrangian relaxation method. To improve the decomposition mechanism, a multi-Pareto cut version is applied in the proposed algorithm. In our numerical tests a modification of the well-known CAB data set is used with different levels of parameters uncertainty. The results demonstrate the capability of the proposed model to design a robust network. We also verify the accuracy of the sample average approximation method. Finally, the results of the proposed algorithm for different instances were compared to other solution approaches which confirm the efficiency of the proposed solution method.
\end{abstract}

\begin{keyword}
\texttt{Hub-networks, Stochastic risk factors, Sample average approximation, Relax-and-decomposition heuristic}
\MSC[2010] 00-01\sep  99-00
\end{keyword}

\end{frontmatter}


\section{Introduction}\label{se1}

A hub is a switching node which reduces transportation cost among several nodes by reducing number of connections. Due to the underlying uncertainty in real world problems, (i.e. volcanic eruption and weather conditions), models dealing with uncertainty in hub network design has recently attracted attention. Uncertainty can be classified into two categories: stochastic uncertainty and the measurement uncertainty of robust optimization. In stochastic optimization, the value of parameters can be modeled by a probability distribution while a probability distribution cannot be used to model the changing parameters of robust optimization (Contreras et al., 2011a). Because of its applicability to the hub location problem, investigation of robust optimization methods could have greater applicability and we intend to explore its potential in this paper. 

Another challenging issue among researchers is the development of an efficient solution procedure capable to obtaining high quality solution in a reasonable time. Only a few studies have addressed this issue. In this regard our contribution will include the use of an augmented solution procedure based on a modified Lagrangian relaxation (LR) method and this is shown to be efficient in solving large scale problems. 

To the best of authors' knowledge, a systematic robust hub location model that directly considers the affect network risk factors have has not been considered in the literature. In summary, a robust HLP is analyzed using a min-max regret approach to deal with network risk as one of objective functions. The solution procedure will be based on the coupling of a Benders Decomposition (BD) algorithm with LR method in order to handle large scale instances. The main contribution of this paper is twofold. First, we use a probabilistic objective function in a min-max regret form to formulate the so called $p$-model. It seems that defining a robust model utilizing the above features can ensure reliability and result in a sustainable hub network design. Secondly, the sample average approximation (SAA) scheme is applied to generate scenarios and then a relax-and-decompose heuristic is implemented as a solution procedure.

The rest of this paper is organized as follows: section \ref{se2} consists of a review on the literature of uncertainty, heuristic methods based on BD and LR and SAA in the hub location problems. A non-linear robust HLP in a risky environment is formulated as a $p$-model in section \ref{se3}. In section \ref{se4}, the relax-and-decomposition solution methodology is presented. In section \ref{se5}, computational results are demonstrated and finally our conclusions and future research are discussed in the last section.

\section{Literature review}\label{se2}

In this section, previous studies on related aspects of current research are briefly reviewed. That is, uncertainty in the hub location, heuristic methods based on BD and LR and finally SAA in the network design are surveyed.

\subsection{Uncertainty in hub location problems}

Real world applications confirm that many parameters related to the network design are uncertain, so classifying location problems in deterministic and non-deterministic categories is reasonable. In this regard, a brief review of risk in HLP and stochastic HLP is presented as follow:

\subsubsection{Risk in hub location problems}

Some events such as the European ash-cloud and the piracy attacks offshore of Somalia confirm that transshipment may be interrupted because of unpredictable and uncontrolled risk factors. Thus a definitive categorization of risk is difficult to make but recently a review paper by Heckmann et al. (2015) has covered different approaches to supply chain risks. They concluded that it is difficult to present a general definition of risk, due to the existence of diverse viewpoints. Chen et al. (2011) considered some factors such as lightning, earthquakes and sandstorms as environmental risk and then extended an analytic network process (ANP) for an international airport project under these natural disasters. Pishvaee et al. (2014) investigated a probabilistic model for a sustainable medical supply chain network where social, economic and environmental criteria are considered as objective functions. Here an accelerated BD algorithm was applied to solve the model.

There are many studies that assume disruptions as a risk and try to design a sustainable supply chain network (see Atoei et al. 2013; Snyder et al. 2006; Goh et al. 2007). Jabbarzadeh et al. (2012) presented a scenario based supply chain model to deal with the risk of disruption. Also a LR and a genetic algorithm (GA) are used to solve the model. Garcia-Herreros et al. (2014) studied a resilient supply chain where the risk of disruption was defined as the fraction of unavailable time for distribution centres. The model was solved by a multi-cut BD method. In the literature there are only a limited number of studies of the risk in hub networks.

\subsubsection{Stochastic \rm{\&} robust hub location problems}

Vasconcelos et al. (2011) used a decentralized management approach to locate hub facilities. The model was formulated assuming a stochastic link cost. Yang (2009) extended a stochastic HLP where the demand value is a random variable and this varied within three scenarios. Zhai et al. (2012) improved a $p$-model to the minimize network risk. The demand parameter was assumed to be a random variable that followed a probability distribution function. A branch and bound (B\&B) algorithm is applied as a solution method. A two stage stochastic HLP was developed for the air network in Iran by Adibi and Razmi (2015). Here, the demand and transportation cost are assumed to be uncertain parameters in their model. de Camargo et al. (2011) modeled a stochastic single allocation HLP in a queueing system. A coupling BD method with Outer Approximations (OA) was used as a solution procedure. Contreras et al. (2011a) constructed a scenario based HLP assuming stochastic demands and costs. A SAA scheme was employed as a scenario generation method and the model was solved by the BD algorithm.

Furthermore, Marianov and Serra (2003) modeled airports as hubs in a M/D/c queuing system and a chance constraint was defined for the number of hub facilities in the queueing system. The model was solved using a tabu search (TS) algorithm. Sim et al. (2009) developed a stochastic $p$-hub center problem to minimize travel time. Travel time was assumed to be a random variable and was considered in a chance constraint form. Hult et al. (2014) proposed a stochastic uncapacitated single allocation $p$-hub center problem. A chance constraint was defined to model the travel time. 

Boukani et al. (2014) published a robust capacitated HLP while fixed cost and hub capacity are assumed as uncertain parameters. The regret measure is used as a robust approach in their model. Shahabi and Unnikrishnan (2014) introduced a robust capacitated HLP where demand is uncertain and is affected by some risk factors such as changes in economy, emergence of new companies as competitors and changes in policy. Alumur et al. (2012) constructed a scenario based HLP to minimize the maximum regret where demand and cost parameters are based on scenarios.

\subsection{Relax-and-decomposition heuristics in hub location problems}

A brief review of previous studies demonstrates that solving the HLP is a challenging problem that has only been addressed during last decade. For this reason, many algorithms have been extended to solve this problem and recently some heuristics based on exact solution methods have been applied. In the following, a brief review of heuristic methods based on LR and BD algorithms is presented. 

Lagrangian relaxation is a methodology to find a bound on the optimal value of a problem for large scale instances in which their optimal solution cannot be found exactly in a reasonable time. The time consuming nature of exact solution approaches for the HLP, have motivated many authors to apply LR (see Lee et al., 1996; Zheng et al., 2014). Contreras et al. (2009) extended LR to obtain tight bounds for the capacitated HLP. We will show that the proposed method improves the best known solution for large size instances.

Due to the time limit, memory issues and the slow convergence of classical LR methods, researchers have recently tried to improve classical versions of solution methods. In this regard, some of major studies are mentioned in the following. Ishfaq and Sox (2011) introduced a hub location-allocation in intermodal logistic networks. To solve the model, a tabu search meta-heuristic and a LR method were applied and results were compared. Lin et al. (2012) presented a capacitated $p$-hub median problem with integral constraints. A GA was coded to solve the model. Also a LR method is applied to show the efficiency of the GA.

LabbÃ© and Yaman (2008) applied B\& B algorithm and LR method to solve the HLP in a star-star network. MariÂ´n (2005) investigated LR method and branching method as a solution procedure to solve the uncapacitated Euclidean HLP. Snyder et al. (2007) proposed a scenario based supply chain model considering risk pooling. A LR method is applied to solve the model. Also a variable fixing approach and B\& B method are employed to close the gap when LR method cannot reach the optimal solution. Yaman (2008) proposed a star $p$-hub median problem with modular arc capacities. A heuristic algorithm based on LR method and local search was applied as a solution method. An et al. (2015) developed a reliable HLP in which the failure of hubs, backup hubs and alternative links are considered. A LR method and a B\& B were used to solve the model. Contreras et al. (2011) presented a multi-period HLP and applied a B\& B algorithm as a solution method which is improved by means of LR method.

Gelareh et al. (2010) presented a shipping liner network hub design model in a competitive environment. An accelerated Lagrangian method embedded in a primal heuristic was considered to solve the model. Karimi and Setak (2014) surveyed a multiple allocation HLP under the incomplete hub location-routing network design. Some lower bounds were obtained by LR and a set of valid inequalities. Lu and Ting (2013) presented LR for the capacitated single allocation $p$-hub median problem. Cacchiani et al. (2013) provided a LR based heuristic for a train-unit assignment problem where the LR method was applied in the model and the relaxed problem was decomposed into a set of assignment problems. Weng and Xu (2014) improved hub routing problem of merged tasks in a collaborative transportation problem. Two heuristic approaches based on LR and BD methods were provided. He et al. (2015) developed an improved MIP heuristic for the intermodal HLP. Their heuristic combines B\&B, Lagrangian relaxation and a linear programming relaxation (LP).

We review previous research in hub location problems that use the BD algorithm. The BD method is one of the exact solution methods used to solve a model when the model contains complicating variables. The BD method has been employed to solve the HLPs (de Camargo et al., 2008; Gelareh and Nickel, 2008; de Camargo et al., 2009a; de Camargo et al., 2009b; Contreras et al., 2011b; Contreras et al., 2012; de SÙ„ et al., 2013), as well as many other problems (Montemanni and Gambardella, 2005; Gendron, 2011). 

Recently authors focus on extended solution methods of decomposition based algorithms to appraise their effectiveness. In some cases of the BD method, it is observed that for the solution of problems (master and sub-problems) the quality of cuts are one of the main issues which affect to convergence rate (Saharidis et al. 2010). Magnanti and Wong (1981) faced, in their seminal work, this same issue in network design problems where they introduced a multi-optimal cuts procedure. A set of cuts, namely Pareto-optimal cuts, are generated in each iteration, where Pareto-cut refer to a cut which is not dominated by other cuts. The results obtained confirmed that these additional cuts improve the convergence of the algorithm significantly. Later, some authors such as Papadakos (2008), Sherali and Lunday (2013) and Tang et al. (2013) employed this accelerating technique in their problems. Also, de Camargo et al. (2013) applied a modified BD to a number of routing hub location problem in which, the Pareto-optimal cut in a B\& B framework was used to solve the master problem. de SÙ„ et al. (2013) produced an accelerated BD method for the tree of the HLP. A cut selection scheme based on Pareto-optimal cuts was proposed for additional efficiency. O'Kelly et al (2014) deployed an improved BD for HLP with price-sensitive demands. de SÙ„ et al. (2015) produced a Benders-branch-and-cut algorithm and several heuristic algorithms including a variable neighborhood descent method, a greedy randomized adaptive search procedure (GRASP) and an adaptive large neighborhood search to solve the $q$-line HLP. Gelareh et al. (2015) presented a multi-period HLP where hub facilities can be opened and closed over different period times. They used an extended metaheuristic solution algorithm with an improved BD method to solve their model. We see that, most of the previous studies related to BD algorithm have focused on improving branching techniques. 

There are some studies that employed BD and LR methods separately, but it is expected that an aggregated algorithm may be more helpful to find a solution efficiently. Gendron (2011) considered a multi-commodity capacitated fixed charge network design problem. Three classes of methods where considered include a cutting-plane method, a BD algorithm and a LR approach which were used to solve the model. Ketabchi and Behboodi-Kahoo (2015) proposed a quadratic L-shaped method based on the projection and the augmented Lagrangian methods for stochastic linear programming. Gelareh and Nickel (2011) investigated an uncapacitated multiple allocation hub location problem tailored for urban transport and liner shipping networks. An accelerated BD and a local search heuristic are employed as solution methodology.

\subsection{SAA in the flow network design}

The sample average approximation technique is a procedure based on Monte Carlo simulation to approximate the stochastic problems with large numbers of scenarios. Calculating the expected value of the objective function is frequently a challenge in the stochastic problems. SAA generates a set of samples of the stochastic parameters to approximate the expected objective value (Kleywegt et al. 2002). 

Santoso et al. (2005) investigated a stochastic supply chain problem and a solution approach based on SAA and BD methods. SchÃ¼tz et al. (2009) introduced a stochastic supply chain in Norwegian meat industry. They used a dual decomposition method along with SAA to solve their model. Bidhandi and Yusuff (2011) developed an integrated stochastic supply chain model as a two stage stochastic program. The model is solved using the SAA couple with BD method. Wang et al. (2013) provided liner ship fleet deployment problem as a joint chance constrained programming model. A SAA scheme was utilized as solution method and the model was solved by SAA and BD methods. 

Table \ref{ta1} is a summary of major previous researches which highlights contribution of the current study. 

\begin{sidewaystable}
\renewcommand{\arraystretch}{.7}
\tabcolsep=3pt
\scalebox{.5}{
\begin{tabular}{ccp{3cm}p{2cm}p{2cm}p{2cm}p{3cm}p{3cm}p{2cm}p{2cm}cp{2cm}p{2cm}} \hline 
\multirow{3}{*}{Reference} & \multirow{3}{*}{Year} & \multirow{3}{*}{Model} & 
\multirow{3}{*}{\rotatebox{90}{no. objectives}}&  \multirow{3}{*}{\rotatebox{90}{MODMA\footnote{Multi-objective decision making approach}}} & \multirow{3}{*}{\rotatebox{90}{Risk type}} & \multirow{3}{*}{\rotatebox{90}{Risk factors}} & \multicolumn{2}{c}{Uncertainty} & {{Sampling method}} & \multicolumn{3}{c}{Solution approach} \\ \cline{8-9}\cline{11-13}
 &  &  &  &  &  &  & Parameter & Approach &  & \multicolumn{2}{c}{Purely solution method} & Hybrid method \\ \hline 
 &  &  &  &  &  &  &  &  &  & {{Exact soln method}} & {{Metaheuristic algorithm}} &  \\ \hline 
Jabbarzadeh et al. & 2012 & Supply chain network & Single & - & Expected value & - & Cost & Scenario-based & - & Classical LR & GA & - \\ \hline 
Atoei et al. & 2013 & Supply chain network & Multiple & $\varepsilon $-constraint method & Expected value & Environmental, macro and supply risks & Reliability/Capacity & Scenario-based & - & - & NSGA-II & - \\ \hline 
Garcia-Herreros et al. & 2014 & Supply chain network & Single & - & Expected value & Operational risks & The availability of DC & Scenario-based & - & Multi-cut BD & - & - \\ \hline 
Pishvaee et al. & 2014 & Supply chain network & Multiple & WSM\footnote{ $ Weighted sum method$ } & possibility & Environmental, macro and security risks & - & - & - & Accelerated BD & - & - \\ \hline 
Marianov and Serra & 2003 & HLP & Single & - & - & - & Number of plane & Probability & - & - & TS & - \\ \hline 
Snyder et al. & 2007 & Location model & Single & - & - & - & Cost & Scenario-based & - & Classical LR & - & - \\ \hline 
Camargo et al. & 2008 & HLP & Single & - & - & - & - & - & - & Classical BD & - & - \\ \hline 
Sim et al. & 2009 & $p$-hub center problem & Single & - & - & - & Time & Chance constraint & - & - & - & Radial heuristic \& Teitz--Bart heuristic \\ \hline 
Yang & 2009 & HLP & Single & - & - & - & Demand & Scenario-based & - & - & - & Heuristic methods \\ \hline 
Camargo et al. & 2011 & HLP & Single & - & - & - & Demand/ Cost & Probability &  & - & - & BD \& OA \\ \hline 
Contreras et al. & 2011a & HLP & Single & - & - & - & Demand/ Cost & Stochastic & SAA & Multi-cut BD & - & - \\ \hline 
Alumur et al. & 2012 & HLP & Single & - & - & - & Demand/ Cost & Scenario-based & - & - & - & - \\ \hline 
Zhai et al. & 2012 & HLP & Single & - & Probability function & Macro risks & Demand & Chance constraint & - & B\&B & - & - \\ \hline 
He et al. & 2015 & Intermodal HLP & Single & - & - & - & - & - & - & - & - & Heuristic combining B\&B, LR, LP \\ \hline 
Proposed model & 2016 & HLP & Multiple & WSM & Regret/ Expected value & Environmental, operational, macro, security and supply risks & Cost/Risk factors & Scenario-based & SAA & - & - & Relax-and-decomposition heuristic \\ \hline 
\end{tabular}}
\caption{A brief review on recent researches of non-linear $p$-robust HLP}\label{ta1}
\end{sidewaystable}

Considering the uncertain nature of real applications (such as equipment breakdowns, ash storm) the use of stochastic variables, makes simulation closer to reality. Also, the memory and time limitations imposed by the use of exact solution algorithms, motivates us to employ a hybrid methodology.

Our survey of the literature reveals that robust approaches to the study of HLPs have received less attention especially when considering network risks in a stochastic environment. The potential to improve relax-and-decomposition heuristics could be an attractive approach in solving large scale HLPs. In this regard, a non-linear robust $p$-model HLP is studied in a risky environment where we employ augmented chance constraint in a min-max regret form to consider network risk as one of the objective functions. Then SAA scheme is applied to select scenarios. One of our main goals is to make use of an augmented BD algorithm within a LR method to improve the solution times and convergence. The complicated constraints are relaxed by LR method and then BD method is used to solve the relaxed MIP. The details of the improved algorithm are described in section \ref{se4} and the potential of this methodology as an effective method is confirmed by numerical results. In summary, we mention our main innovations as follow:

\begin{itemize}
\item  Formulating a $p$-model approach in a $\min-\max$ regret form. 
\item  Directly considering network risks as stochastic factors.
\item  Applying SAA scheme to generate scenarios.
\item  Improving a relax-and-decomposition heuristic as a solution procedure.
\end{itemize}

It seems that the proposed model and solution procedure are capable of designing a sustainable hub network and solving large scale problems, respectively.

\subsection{Motivation for the study}

As mentioned earlier, some events such as extreme weather conditions, changes in economy and link failures are examples of unpredictable nature of real-world situations. On the other hand, there are many real-world events that show the drawback of a robust model in that they can led to loss of human life and financial losses. The occurrence of these events motivates us to propose our $p$-robust HLP model. Followings are examples of these events:

\begin{itemize}
\item  The Valigonda rail disaster occurred in the India on 29 October 2005. A flash flood swept away a small rail bridge and a train traveling on it derailed at the broken section of the line. Resulting, at least 114 and 200 people are killed and injured, respectively.
\end{itemize}

(\url{https://en.m.wikipedia.org/wiki/Valigonda\_train\_wreck})

\begin{itemize}
\item  Following the volcanic ash ejected during the 2010 eruptions of EyjafjallajÙkull in Iceland, many airports including two international hub airports in Frankfurt and Paris, were closed. Over 95,000 flights had been cancelled which result in the largest air-traffic shut-down since World War II. The international air transport association stated that the airline industry worldwide would lose around \$1.7 billion during the disruption.
\end{itemize}
(\url{https://en.wikipedia.org/wiki/Air\_travel\_disruption\_after\_the\_2010\_Eyjafjallaj
\%C3\%B6kull\_eruption})

\begin{itemize}
\item
Damaging major roads and closing some airports such as Louis Armstrong New Orleans International Airport were some effects of Hurricane Katrina in New Orleans in 2010. The disaster had major implications for a large segment of the population, economy, and transportation of the entire United States.
\end{itemize}

(\url{https://en.wikipedia.org/wiki/Effects\_of\_Hurricane\_Katrina\_in\_ew\_Orleans})

It is noteworthy that containerized cargo on vessels are employed as a hub-and-spoke transportation system for up to 90 percent of global trade volume (Gelareh and Nickel 2011). The liner shipping industries lease container terminals for a certain time horizon with the objective to minimize the waiting times of the vessels. In this particular case, the use of network design for a relatively short period of time makes sense because of the flexibility it provides in relocating their calling ports. Depending on their future market share, a multi-period modeling framework provides a better decision framework than the static counterpart.

Previous real events especially on the maritime transportation networks motivated us to consider a multi-period network design problem in presence of risk factors.

\section{The proposed model for the $p$-robust hub location problem }\label{se3}
\subsection{Model description}

Due to the lack of information about effected parameters in real applications, it is difficult to manipulate a probability distribution function. On the other hand, lack of enough data makes the expected cost measure useless. In this regard, a min-max regret can be used as robust measure to consider network risk under uncertainty which can then formulated in a $p$-model form. The proposed model minimizes the network cost and the maximum regret deals with network risk. We intend to design an inter-hub network where we have only a few similar works such as Contreras et al. (2011), Zhai et al. (2012) and Gelareh et al. (2015) which study this kind of hub networks. It is worth noting that the probability of the network risk free scores to be less than a threshold value is considered as an objective function. In the following model is described.

\noindent\scalebox{.9}{\begin{tabular}{lp{12cm}} 
\multicolumn{2}{l}{\textit{Indices}}\\
\hline 
$i,j$ & Potential hub nodes ($i,j\in H)$).\\ 
$s$ & Scenarios ($s,s'\in S$).\\ 
$t$ & Period times ($t\in T$).\\
\hline 
%
\multicolumn{2}{l}{\textit{Parameters}}
\\
\hline 
$f_{i}^{t}$& Fixed cost for locating a hub facility at node~$i\in H$~{in period time} $t\in T$.\\ 
$c_{ij}^{st}$& $\begin{array}{l}{\rm Stochastic~unit~ transportation~ cost~ from~ hub} ~i\in H~{\rm to}~j \in H~{\rm in~ period}\\{\rm time}~ t\in T{\rm based~on~ scenario}~s\in S.\end{array}$\\ 
$w_{i}^{st} $& $\begin{array}{l}{\rm Stochastic ~flow~ to~ be~ routed~ through~ hub ~node}~i\in H~{\rm  in~ period~ time}\\t\in T~{\rm based~on}{\rm scenario}~s\in S.\end{array}$\\ 
$\pi_{1i}^{st} $ & $\begin{array}{l}{\rm Stochastic~risk~ free~ score~factor~ for~ potential~ hub~ node}~i\in H~{\rm in~ period}\\{\rm time}~t\in T{\rm based~on~ scenario}~s\in S.\end{array}$ \\
$\pi_{2ij}^{st} $& $\begin{array}{l}\text{Stochastic risk free score factor for a hub link between hubs}~i\in H\\ \text{and}~j\in H \text{in period time}~t\in T~\text{based on scenario}~s\in S.\end{array}$\\
$\pi_{0}^{1s} $ & Threshold value of sustainability for potential hub nodes based on scenario~$s\in S$.\\
$\pi_{0}^{2s} $ & Threshold value of sustainability for hub links based on scenario $s\in S$.\\
$p^{s} $& Probability of occurrence of scenario $s\in S (\sum_{s\in S}p^{s} =1)$.\\ \hline 
\multicolumn{2}{l}{\textit{Decision variables}}
\\
 \hline 
$Z_{i}^{t} $ & Is equal to 1, if hub facility is located in nod~ $i\in H$ in period time~ $t\in T$. \\ 
$X_{ij}^{st} $ &\text{Is equal to 1, if an inter-hub connection between hubs}~ $i\in H$
{ and}~ $j\in H$~{is constructed in period time it}~ $t\in T$~{ based on} {scen}~ $s\in S$.\\ \hline 
\end{tabular}}

Note that, changing the data underlying the problem over a long-time period convinced us to consider the multi-period nature of this problem. To disregard this aspect of the problem may result in enduring additional costs for the decision maker. Hence, hub facilities are allowed to be opened and closed at different time periods to enhance a flexible in the transportation system. Let us also denote the total set of nodes by G. Some seasonal criteria such as demand within season, seasonal business targets, seasonal passengers demand (when dealing with tourist areas) and specialization of ships (special ships like tankers, grain carriers, barges, mineral carriers, bulk carriers and container ships) lead to limitations, for the shipping industries, in having certain selections of vessels for a short horizons of time. For this reason, we assume that the allocated origin $o(t)\in G$ and destination $d(t)\in G$ to a container terminals are known for a certain period of time $t\in T$. Accordingly, the stochastic unit transportation cost in period time $t\in T$, based on scenario $s\in S$, is defined as $c_{ij}^{st} =c_{o(t)i}^{'st} +\tau c_{ij}^{'st} +c_{jd(t)}^{'st} $ where $c'$ and $\tau $ are referred to transportation costs and transportation discount factor among hub links along the path  $\left(o(t), i, j, d(t)\right)$, respectively. Also, if a hub facility has been opened at a specific period time, it will serve just for that period time. This is considered as another assumption in the proposed model. The proposed model is formulated as a non-linear two-stage stochastic program based on different scenarios. The first stage decisions are location of hub nodes and then in the second stage the hub links are assigned for each scenario $s\in S$.

The steps of the modeling are organized as follows: 
\begin{itemize}
\item[Step 1:]
Formulating the problem in a two stage $p$-robust stochastic HLP.
\item[Step 2:]
Linearizing the non-linear program model.
\item[Step 3:]
Using the WSM method to defining the equivalent integrated single objective model.
\end{itemize}
Firstly, stochastic risk free score factors corresponding to different types of factors which should be integrated into aggregated indices as input data. The value of these factors are calculated based on the score values which are considered for them. The calculated risk free score factors are used to define the risk objective function. The maximum difference between network risks and the best objective function value of each scenario is minimizes as follows:

\begin{equation}\label{eq1}
\mathop{\min\max\partial_s}_{s\in S} 
\end{equation}

In which:

\begin{align}\label{eq2}
\partial_s=&\mathbf{pr}\left(\sum_{i\in H}\sum_{t\in T}\pi^{st}_{1i}(1-Z^{t-1}_i)Z^t_i\leq \mathbf{\pi^{1s}_0}\right)\mathbf{ pr}\left(\sum_{i\in H}\sum_{j\in H}\sum_{t\in T}\pi^{st}_{2ij}X^{st}_{ij}\leq \mathbf{\pi^{2s}_0}\right)\nonumber\\
&\qquad\qquad\qquad\qquad\qquad\qquad\qquad\qquad\qquad\qquad -\Psi^{**}_s\qquad s\in S
\end{align} 

where $\Psi^{**}_s$ is achieved by the following model and $\partial_s$ is the regret value related to the risk objective for each scenario $s\in S$. It is noteworthy that as each node and connection risks are independent, their probabilities can be multiplied in order to be jointly considered.  

\begin{align}
&\min \mathbf{ pr}\left(\sum_{i\in H}\sum_{t\in T}\pi^{st}_{1i}(1-Z^{t-1}_i)Z^t_i\leq \mathbf{ \pi^{1s}_{0}}\right)\mathbf{ pr}\left(\sum_{i\in H}\sum_{j\in H}\sum_{t\in T}\pi^{st}_{2ij}X^{st}_{ij}\leq \mathbf{ \pi ^{2s}_{0}}\right)\label{eq3}
\end{align}
Subject to:
\begin{align}
&\sum_{i\in H}{\sum_{j\in H}{X^{st}_{ij}}}=1&&s\in S, t\in T\label{eq4} \\
&\sum_{j\in H}X^{st}_{ij}+\sum_{j\in H\backslash\{i\}}X^{st}_{ji}\leq Z^t_i&&i\in H, s\in S, t\in T\label{eq5}\\
&X^{k\ s}_{ij}\geq 0, Z_i\in {\mathbb{ B}}^{|H|}&&i,j\in H,s\in S,t\in T\label{eq6}
\end{align} 

The probability that the network total risk free scores are not less than a threshold value is minimized as an objective function in equation \eqref{eq3}. Constraints \eqref{eq4} guarantee that one hub path should be allocated based on scenario $s\in S$ in each time period $t\in T$ in order that all demand is fully routed through the network. Constraints \eqref{eq5} ensures that if a hub facility is located at node $i\in H$, then both collection and distribution can occur. Finally, constraints \eqref{eq6} define the decision variables. To simplify the complexity of the model we only consider the variables $X^{ks}_{ij}$ to be continuous. On the other hand the values $Z^t_i$ are binary and the coefficient matrix associated with constraints \eqref{eq4}-\eqref{eq5} is totally unimodular, so the continuous relaxation of problem \eqref{eq3}-\eqref{eq6} will always have a binary optimal solution.

Now, the non-linear multi-objective $p$-robust hub location problem, considering network risks (PRH-R), can be defined as follows:

\begin{align}
&{\rm Min}\mathop{\max\partial_s}_{s\in S}\label{eq7}\\
&{\rm Min}~\Omega =\sum_{t\in T}\sum_{i\in H}f^t_i(1-Z^{t-1}_i)Z^t_i+\sum_{t\in T}\sum_{i\in H}\sum_{j\in H}\sum_{s\in S}p^sw^{st}_i c^{st}_{ij} X^{st}_{ij}\label{eq8}
\end{align}
Subject to:\eqref{eq4}-\eqref{eq5}
\[
X^{s t}_{ij}\geq 0, Z_i\in \mathbb{B}^{|H|} \qquad i, j\in H, s\in S, t\in T\nonumber
\]

The objective function is defined by equations \eqref{eq7}-\eqref{eq8}. As earlier mentioned, the first one minimizes the maximum differences between network risks and the best objective function value of each scenario. The second one minimizes expected network costs in different time period and scenarios where the recovery gain generated by closing a hub facility in each time period is considered in its setup cost. 

\subsubsection{Linearization of the proposed model}
As stated before there are some non-linear terms in the proposed model, however some transformations should be done to acquire a standard MIP form. We assume a multivariate normal distribution for the $|H|\times |T|-$dimensional random vector $\pi_{1} =\left(\pi_{11}^{1}(s),\pi_{11}^{2} (s),\ldots,\pi_{1|H|}^{|T|} (s)\right)^{T} .$ Therefore, an $|H|\times |T|\times s$ matrix $A$ and a vector $\mu =\left(\mu_{1}^{1} ,\mu_{1}^{2} ,\ldots,\mu_{|i|}^{|t|} \right)\in \mathbb{ R}^{|i|\times |t|} $exist in which $\pi_{1} =AS'+\mu $ where $S'$ is an $s$-dimensional independent random vector with a standard normal distribution (Zhai et al. 2012). Hence, the covariance matrix of $\pi_{1} $ is $\Sigma =AA^{T} $. By defining $\wp (Z,S)=\sum _{i\in H}\sum _{t\in T}\pi _{1i}^{st} \left(1-Z_{i}^{t-1} \right)  Z_{i}^{t} $, we can represent $\wp (Z,S)$ as $\wp (Z,S)=\pi _{1}^{T} Z=S'^{T} \left(A^{T} Z\right)+\mu ^{T} Z$ where $Z=\left(Z_{1}^{1}, Z_{1}^{2}, \ldots, Z_{|H|}^{|T|} \right).$ Now, the part of the risk objective function related to hub facilities can be reformulated using $E \left(\wp (Z,S)\right)=\mu^{T} Z=\sum _{i\in H}\sum _{t\in T}p^{s} \pi _{1i}^{st} \left(1-Z_{i}^{t-1} \right)  Z_{i}^{t} $ and ${ Var}\left(\wp (Z,S)\right)=\|A^{T} Z\|^{2} =\sigma _{{\pi}_{1}^{st} }^{2} $ as:

\begin{align}
pr\left(\wp (Z,S)\leq \pi _{0}^{1s} \right)&=pr\left(\frac{\wp (Z,S)-E \left(\wp (Z,S)\right)}{\sqrt{{\rm Var}\left(\wp (Z,S)\right)} } \leq \frac{\pi _{0}^{1s}-E \left(\wp (Z,S)\right)}{\sqrt{{\rm Var}\left(\wp (Z,S)\right)} } \right)\label{eq9}\\
= &pr\left(\frac{\wp (Z,S)-\sum _{i\in H}\sum _{t\in T}p^{s} \pi _{1}^{st} \left(1-Z_{i}^{t-1} \right)  Z_{i}^{t} }{\sigma _{{\pi }_{1}^{st} } } \leq\right.\nonumber\\
&\qquad\qquad\qquad\quad \left.\frac{\pi _{0}^{1s} -\sum _{i\in H}\sum _{t\in T}p^{s} \pi _{1i}^{st} \left(1-Z_{i}^{t-1} \right)  Z_{i}^{t} }{\sigma _{{\pi }_{1}^{st} } } \right)\nonumber\\
=&\phi _{1} \left(\frac{\pi _{0}^{1s} -\sum _{i\in H}\sum _{t\in T}p^{s} \pi _{1i}^{st} \left(1-Z_{i}^{t-1} \right)  Z_{i}^{t} }{\sigma _{{\pi }_{1}^{st} } } \right)\nonumber
\end{align}
Similarly, in the same way the second part of the risk objective function can be reformulated as:
\begin{equation}\label{eq10}
pr\left(\sum _{i\in H}\sum _{j\in H}\sum _{t\in T}\pi _{2ij}^{st}    X_{ij}^{st} \le \pi _{0}^{2s} \right)=\phi _{2} \left(\frac{\pi _{0}^{2s} -\sum _{i\in H}\sum _{j\in H}\sum _{t\in T}p^{s} \pi _{2ij}^{st}    X_{ij}^{st} }{\sigma _{{\pi }_{2}^{st} } } \right)
\end{equation}
Now, equation \eqref{eq2} can be changed to the following equation.

\begin{align}\label{eq11}
{\partial }_s=&{\phi }_1\left(\frac{{\mathbf{ \pi }}^{\mathbf{ 1}\mathbf{ s}}_{\mathbf{ 0}}\mathbf{ -}\sum_{i\in H}{\sum_{t\in T}{p^s{\pi }^{st}_{1i}\left(1-Z^{t-1}_i\right)Z^t_i}}}{{\sigma }_{{\pi }^{st}_{1}}}\right)\nonumber\\
&\times\phi_2\left(\frac{{\mathbf{ \pi }}^{\mathbf{ 2}\mathbf{ s}}_{\mathbf{ 0}}\mathbf{ -}\sum_{i\in H}{\sum_{j\in H}{\sum_{\mathbf{ t}\in T}{p^s{\pi }^{st}_{2ij}X^{st}_{ij}}}}}{{\sigma }_{{\pi }^{st}_{2}}}\right)-{\Psi }^{**}_s\qquad s\in S
\end{align}

If the amount of probabilities is more than 0.5 then probabilities will be increased by increasing of their standard values, so equation \eqref{eq11} can be substituted into equation \eqref{eq12}.

\begin{align}\label{eq12}
{\partial }_s=&\frac{{\mathbf{ \pi }}^{\mathbf{ 1}\mathbf{ s}}_{\mathbf{ 0}}\mathbf{ -}\sum_{i\in H}{\sum_{t\in T}{p^s{\pi }^{st}_{1i}(1-Z^{t-1}_i)Z^t_i}}}{{\sigma }_{{\pi }^{st}_{1}}}\nonumber\\
&\mathbf{ \times }\frac{{\mathbf{ \pi }}^{\mathbf{ 2}\mathbf{ s}}_{\mathbf{ 0}}\mathbf{ -}\sum_{i\in H}{\sum_{j\in H}{\sum_{\mathbf{ t}\in T}{p^s{\pi }^{st}_{2ij}X^{st}_{ij}}}}}{{\sigma }_{{\pi }^{st}_{2}}}-{\Psi }^*_s\qquad s\in S
\end{align}

where ${\Psi }^*_s$ is updated to the optimal value of linearized objective function considering constraints \eqref{eq4}-\eqref{eq6}. In spite of above linearization, several nonlinear terms are appeared in the proposed model consisting of products of binary variables and the $\min -\max$ form. Therefore we replace $V_{i}^{t,t-1} =Z_{i}^{t-1} Z_{i}^{t}$, $L_{ij}^{st} =Z_{i}^{t} X_{ij}^{st}$ and $Q_{ij}^{st, t-1} =V_{i}^{t ,t-1} X_{ij}^{st}$ in the model and add corresponding auxiliary constraints \eqref{eq16}-\eqref{eq21}. Also axillary variable $\gamma $ is defined to linearize the min-max model as follows (Glover and Woolsey, 1974):

\begin{align}
&{\rm Min}\ \Omega =\sum_{t\in T}{\sum_{i\in H}{f^t_i\left(Z^t_i-V^{t,\ t-1}_i\right)}}+\sum_{t\in T}\sum_{i\in H}{\sum_{j\in H}{\sum_{s\in S}{p^sw^{st}_i c^{st}_{ij}X^{st}_{ij}}}}\label{eq13}\\
&{\rm Min}~\gamma \label{eq14}
\end{align}
subject  to: \eqref{eq4}-\eqref{eq5}
\begin{align}
\gamma \geq & \frac{{\mathbf{ \pi }}^{\mathbf{ 1}\mathbf{ s}}_{\mathbf{ 0}}{\mathbf{ \pi }}^{\mathbf{ 2}\mathbf{ s}}_{\mathbf{ 0}}\mathbf{ -}{\mathbf{ \pi }}^{\mathbf{ 1}\mathbf{ s}}_{\mathbf{ 0}}\sum_{i\in H}{\sum_{j\in H}{\sum_{\mathbf{ t}\in T}{p^s{\pi }^{st}_{2ij}X^{st}_{ij}}}}\mathbf{ -}{\mathbf{ \pi }}^{\mathbf{ 2}\mathbf{ s}}_{\mathbf{ 0}}\sum_{i\in H}{\sum_{t\in T}{p^s{\pi }^{st}_{1i}\left(Z^t_i-V^{t,\ t-1}_i\right)}}}{{\sigma }_{{\pi }^{st}_{1}}{\sigma }_{{\pi }^{st}_{2}}}&&\nonumber
\end{align}
\begin{align}
&\qquad+\frac{\sum_{i\in H}{\sum_{j\in H}{\sum_{\mathbf{ t}\in T}{p^s{\pi }^{st}_{2ij}{\pi }^{st}_{1i}(L^{st}_{ij}-Q^{st,t-1}_{ij})}}}}{{\sigma }_{{\pi }^{st}_{1}}{\sigma }_{{\pi }^{st}_{2}}}-{\Psi }^*_s &&s\in S\label{eq15}\\
&V^{t,t-1}_i\geq Z^{t-1}_i+Z^t_i-1&& i\in H, t\in T\label{eq16}\\
&2V^{t,t-1}_i\leq Z^{t-1}_i+Z^t_i&& i\in H, t\in T\label{eq17}\\
&L^{st}_{ij}\geq Z^t_i+X^{st}_{ij}-1&&i,j\in H, s\in S, t\in T\label{eq18}\\
&2L^{st}_{ij}\leq Z^t_i+X^{st}_{ij}&&i,j\in H, s\in S, t\in T\label{eq19}\\
&Q^{st,t-1}_{ij}\geq V^{t,t-1}_i+X^{st}_{ij}-1&&i,j\in H, s\in S, t\in T\label{eq20}\\
&2Q^{st,t-1}_{ij}\leq V^{t,t-1}_i+X^{st}_{ij}&&i,j\in H, s\in S, t\in T\label{eq21}\\
&X^{st}_{ij},\gamma \geq 0, Z^t_i,V^{t,t-1}_i,L^{st}_{ij},Q^{st,t-1}_{ij}\in
\mathbb{ B}^{\left|H\right|}&&i,j\in H, s\in S, t\in T\label{eq22}
\end{align} 
Constraint \eqref{eq16} and \eqref{eq17} ensure that the optimal value of variables $V_{i}^{t,t-1} $ are equal to the product of the variables $Z_{i}^{t-1} $ and $Z_{i}^{t} $. Similarly, this situation pertains to the variables $L_{ij}^{st} $ and $Q_{ij}^{st,\, t-1} $ which is handled by the Constraint \eqref{eq18}-\eqref{eq19} and \eqref{eq20}-\eqref{eq21}, respectively.

\subsubsection{Multi-objective optimization methodology}
In this subsection a linear composite objective function, namely, the weighted sum method, is used to define an equivalent integrated single objective model (Deb 2001). Because of existing different units in the objective functions, a scalarization strategy is needed. In this regard, the ideal and nadir values of the objective functions are given by ${\gamma }^*, {\Omega }^*$ and $\gamma^{\max}, {\Omega }^{\max}$ , respectively. Moreover, symbols${\theta }_n(\sum_n{{\theta }_n}=1$, ${\theta }_n>0), {n=1, 2}$ are denoted the weights of objective functions. Now the composite objective function is given by the following equation.

\begin{align}\label{eq23}
&{\rm Min}\ \omega ={\theta }_1\left(\frac{\gamma -{\gamma }^*}{{\gamma }^{max}-{\gamma }^*}\right)\nonumber\\
&+\frac{\theta_2 }{{\Omega }^{max}-{\Omega }^*}\left({\sum_{t\in T}{\sum_{i\in H}{f^t_i\left(Z^t_i-V^{t,t-1}_i\right)}}+\sum_{t\in T}{\sum_{i\in H}{\sum_{j\in H}{\sum_{s\in S}{p^sw^{st}_i{\ c}^{st}_{ij}X^{st}_{ij}}}}}-{\Omega }^*}\right)
\end{align}
Subject to: \eqref{eq4}- \eqref{eq5}, \eqref{eq15}-\eqref{eq21}
\[
X^{st}_{ij},\gamma \ge 0\ ,Z^t_i,V^{t,t-1}_i,L^{st}_{ij},Q^{st,t-1}_{ij}\in {\mathbb{ B}}^{\left|H\right|}\qquad
i,j\in H, s\in S, t\in T\]

where, $\omega $ is the optimum value of the equivalent single objective function.

\section{Solution methodology}\label{se4}

\subsection{Motivation for proposing the relax-and-decomposition heuristic method}

In our case, a five index formulation of the PRH-R led to a huge number of constraints and variables which results in very high complexity even for small instances. Exact solution methods can be time consuming for such problems, hence our focus will be on extending heuristics based on exact solution methods. We propose a relax-and-decompose method to solve our model. The main idea for our solution algorithm is to relax the complicating constraints by a LR method to obtain tight bounds for the proposed model. The slow convergence of the Lagrangian relaxation applied in large-scale instances motivated us to develop a solution methodology to accelerate the classical LR method. The relaxed MIP still contains complicated variables therefore a decomposition method is crucial for this step. Consequently, we proposed a relax-and-decompose approach applied to the LR method in order to remove complicated constraints and then employ a BD algorithm to the remaining Lagrangian subproblem. Moreover, as the number of complicated variables is much smaller than the continuous ones, the efficiency of the BD method is enhanced, while this observation   does not hold for the original model where the Benders master problem can be difficult to solve because of its large size. In this regard, generating strong cuts using an efficient set of the sub-problem solutions can improve the convergence rate. Pareto-optimal cuts are a well-known tool for this purpose. Also, a multi-generation cut procedure can be applied to enhance the algorithm's efficiency, where adding multi-cut in each iteration speeds up the classical BD algorithm. Thus we propose to use a multi-Pareto optimal cut BD method to solve the relaxed MIP. This strategy is capable of solving large-scale instances in a reasonable time.

 Our model is based on random scenarios. So generating effective scenarios is important if they are to cover the real situational features. In this regard, the SAA scheme is one of many random scenario generation methods which can be used. In the following, first the SAA scheme is employed to select scenarios and then the upper bound is calculated by using the following proposed solution procedure. The major steps for the solution method are given in Figure 1. After using the SAA as a sampling method, some constraints which lead to difficulty in solving the PRH-R are relaxed using LR algorithm. Then a multi-Pareto cut BD algorithm is used as an internal procedure to solve the remaining relaxed problem. Notice that, the inner BD produces a lower bound (${\rm LB_{LR}}$) which is considered as the Lagrangian lower bound (${\rm LB_{LR}}$) within the algorithm. Finally, a subgradient method is applied to update Lagrange multiplier sets.

\begin{figure}[h!]
\centering
\includegraphics[scale=0.5]{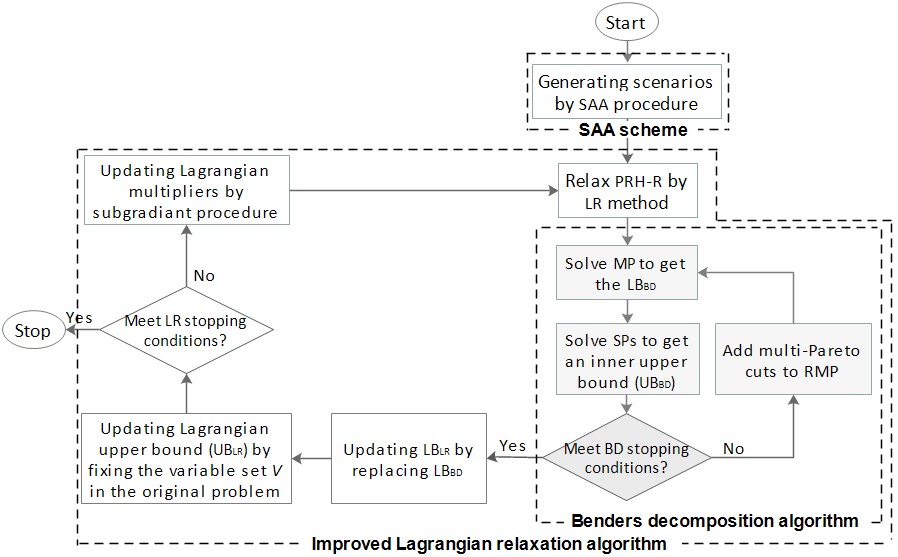}
\caption{The proposed solution heuristic diagram}\label{fig1}
\end{figure}

\subsection{Sample average approximation scheme}

One of difficulties met with chance constraint programming is that the feasible region generally is not convex. For this reason, previous studies have looked at different approaches to employ convex approximation (Ben-Tal and Nemirovski, 2000; Hong et al., 2011) and Monte Carlo simulation (Pagnoncelli et al., 2009) to the approximation of the chance constraint model. Because of the binary variables in our proposed chance constraint based model, the convex approximation approaches cannot be employed. Also, as mentioned earlier, calculating the expectation of scenario based parameters is challenging for this model. To address these difficulties, the SAA scheme is used. A set of random sample $S=\{1, 2, \ldots ,s\}$ is generated and the expected value of the objective function is approximated as follow:

\begin{align}\label{eq24}
&{\rm Min}~\omega ={\theta }_1\left(\frac{\gamma -{\gamma }^*}{{\gamma }^{max}-{\gamma }^*}\right)\nonumber\\
&+\frac{\theta_2}{{\Omega }^{max}-{\Omega }^*} \left( {\sum_{t\in T}{\sum_{i\in H}{f^t_i\left(Z^t_i-V^{t,t-1}_i\right)}}+\frac{1}{|S|}\sum_{t\in T}{\sum_{i\in H}{\sum_{j\in H}{\sum_{s\in S}{w^{st}_i{\ c}^{st}_{ij}X^{st}_{ij}}}}}-{\Omega }^*}\right)\nonumber\\
&
\end{align}

Suppose $\omega ^{S} $ and $(\widehat{Z}^{S} ,\widehat{\gamma }^{S} ,\widehat{V}^{S} ,\widehat{L}^{S} ,\widehat{Q}^{S} )$ is the optimal value and an optimal solution vector of the SAA problem \eqref{eq24}, respectively. Now, for a particular scenario, the problem \eqref{eq24}, \eqref{eq4}-\eqref{eq5}, \eqref{eq15}-\eqref{eq22} can be deterministically solved. 

Under mild regularity conditions, it is expected that by increasing the sample size \textbar \textit{S}\textbar , $\omega ^{S} $ and $(\widehat{Z}^{S} ,\widehat{\gamma }^{S} ,\widehat{V}^{S} ,\widehat{L}^{S} ,\widehat{Q}^{S} )$ we obtain converge to their equivalents with probability one (Kleywegt et al., 2002). An acceptable approach to selecting the sample size $S$  by is taking into account the trade-off between the quality of the solution of the SAA problem \eqref{eq24}, \eqref{eq4}-\eqref{eq5}, \eqref{eq15}-\eqref{eq22} and the computational time. For this, a set of independent sample is generated and the SAA problem \eqref{eq24}, \eqref{eq4}-\eqref{eq5}, \eqref{eq15}-\eqref{eq22} is solved repeatedly instead of solving a large-scale SAA problem. Now, the SAA procedure is described as follow:

\begin{itemize}
\item[\textit{Step 1}:] Consider $M=\{ S_{1} ,\ldots,S_{|M|} \} $ as a set of $M$  replications with sample size $|S|$ i.e. $s_{1}^{m}, \ldots, s_{|S|}^{m} $ for $m\in M$. The SAA problem is solved for $M$ replications.

\begin{align*} 
&{\rm Min}~\omega ={\theta }_1\left(\frac{\gamma -{\gamma }^*}{{\gamma }^{\max}-{\gamma }^*}\right)\\
&+\frac{\theta_2 }{{\Omega }^{max}-{\Omega }^*}\left({\sum_{t\in T}{\sum_{i\in H}{f^t_i\left(Z^t_i-V^{t,t-1}_i\right)}}+\frac{1}{|S|}\sum_{t\in T}\sum_{i\in H}\sum_{j\in H}\sum_{s\in S_m}w^{st}_i c^{st}_{ij}X^{st}_{ij}-{\Omega }^*}\right)
\end{align*} 
Subject to: \eqref{eq4}-\eqref{eq5}, \eqref{eq15}-\eqref{eq22}

Assume, $\omega ^{S_{m} } $ and $(\widehat{Z}^{S_{m} }, \widehat{\gamma }^{S_{m} }, \widehat{V}^{S_{m} }, \widehat{L}^{S_{m} }, \widehat{Q}^{S_{m} } )$ be the optimal objective value and an optimal solution in replication $m\in M$, respectively.

\item[\textit{Step 2: }] We can compute the average of the optimal objective values that where obtained at previous step. It is known that this average provides a lower bound for the optimal value of the original problem \eqref{eq24}, \eqref{eq4}-\eqref{eq5}, \eqref{eq15}-\eqref{eq22}. The average $\mu _{M}^{S} $ and an estimate of corresponding variance $\sigma _{\mu _{M}^{S}}^{2} $are formulated as follow:

\begin{align} 
&\mu _{M}^{S} =\frac{1}{|M|} \sum _{m\in M}\omega ^{S_{m} }\label{eq25}\\
&\sigma _{\mu _{M}^{S}}^{2} =\frac{1}{\left(|M|-1\right)|M|} \sum _{m\in M}\left(\omega ^{S_{m} } -
\mu_{M}^{S} \right)^{2}\label{eq26}
\end{align}

\item[\textit{Step 3}:] An upper bound for the optimal value of the original problem \eqref{eq24}, \eqref{eq4}-\eqref{eq5}, \eqref{eq15}-\eqref{eq22} can be estimated by generating an independent reference sample $S^{'} $. So that the size of reference sample $S^{'} $ is much greater than the sample size of $S$. It is obvious that the SAA problem \eqref{eq24}, \eqref{eq4}-\eqref{eq5}, \eqref{eq15}-\eqref{eq22} can be easily computed for reference sample $S^{'} $ by fixing binary variables of the first stage. We choose one of computed solution in \textit{step 1} to fixed binary variables. Now, the estimated SAA upper bound and its variance can be obtained as follows:

\begin{align}
&\omega_{S'}(\widehat{Z}, \widehat{\gamma }, \widehat{V}, \widehat{L}, \widehat{Q})=\theta_{1} \left(\frac{\gamma -\gamma^{*} }{\gamma^{\max }-\gamma^{*} } \right) \label{eq27}\\
&+\frac{\theta_{2}}{\Omega ^{\max } -\Omega ^{*} }  \left({\sum _{t\in T}\sum _{i\in H}f_{i}^{t} \left(Z_{i}^{t} -V_{i}^{t,t-1} \right)  +\frac{1}{|S^{'} |} \sum _{t\in T}\sum _{i\in H}\sum _{j\in H}
\sum_{s\in S^{'} }w_{i}^{st} c_{ij}^{st} X_{ij}^{st} -\Omega^{*} } \right)\nonumber
\end{align}

\begin{align}\label{eq28}
&\sigma_{S^{'}}^{2} (\widehat{Z},\widehat{\gamma },\widehat{V},\widehat{L},\widehat{Q})=\frac{1}{\left(|S^{'} |-1\right)|S^{'} |} \sum _{s\in S^{'} }\left(\theta _{1} \left(\frac{\gamma -\gamma ^{*} }{\gamma ^{\max } -\gamma ^{*} } \right)\right.\\
&\left. +\theta _{2} \left(\frac{\sum _{t\in T}\sum _{i\in H}f_{i}^{t} \left(Z_{i}^{t} -V_{i}^{t,t-1} \right)  +\sum _{t\in T}\sum _{i\in H}\sum _{j\in H}w_{i}^{st} c_{ij}^{st} X_{ij}^{st}-\Omega ^{*} }{\Omega ^{\max } -\Omega ^{*} } \right)\right.\nonumber\\
&\left.-\omega (\widehat{Z},\widehat{\gamma },\widehat{V},\widehat{L},\widehat{Q})\right)^{2}\nonumber 
\end{align}

\item[Step 4:] estimate the optimality gap and its variance by using the estimated SAA bounds and their variances as follow:

\begin{align}
&gap_{S,M,S^{'}}=\omega _{S^{'} } (\widehat{Z},\widehat{\gamma },\widehat{V},\widehat{L},\widehat{Q})-\mu _{M}^{S} \label{eq29}\\ 
&\sigma _{gap}^{2} =\sigma_{\mu_{M}^{S}}^{2} +\sigma _{S^{'}}^{2} (\widehat{Z},\widehat{\gamma }, \widehat{V}, \widehat{L}, \widehat{Q})\label{eq30}
\end{align}
\end{itemize}

\subsection{Relax-and-decomposition heuristic for PRH-R}

We exploited the LR method, combining the advantages of an acceleration of Benders decomposition algorithm. This subsection is devoted to the description of the augmented BD method integration into LR algorithm to speed up time and convergence. 

\subsubsection{LR framework}

Complicating constraints are a set of constraints which make the problem difficult to solve. We observed that the assignment constraint \eqref{eq4} and risk constraint \eqref{eq15} lead to a computational complexity in the model. These constraints can be thought of as complicating constraints and we note that the resulting relaxed model is easier to solve. Hence, they are relaxed in a Lagrangian manner and added into the objective function by using corresponding Lagrange multipliers $d_{1s}$ and $d_{2st}$ for each $s\in S$ and $t\in T$, respectively. In this way, the Lagrangian relaxation of PRH-R (LRP) is defined as follows:

\begin{align}
&LRP\left[d_1,d_2\right]\nonumber\\
&{\rm Min}\ \zeta ={\theta }_1\left(\frac{\gamma -{\gamma }^*}{{\gamma }^{max}-{\gamma }^*}\right)\nonumber\\
&+\frac{\theta_2 }{{\Omega }^{max}-{\Omega }^*} \left({\sum_{t\in T}{\sum_{i\in H}{f^t_i\left(Z^t_i-V^{t,t-1}_i\right)}}+\sum_{t\in T}{\sum_{i\in H}{\sum_{j\in H}{\sum_{s\in S}{p^sw^{st}_i{\ c}^{st}_{ij}X^{st}_{ij}}}}}-{\Omega }^*}\right)\label{eq31}\\
&+\sum_{s\in S}d_{1s}\left(\left(\frac{{{\mathbf{ \pi }}^{\mathbf{ 1}
\mathbf{ s}}_{\mathbf{ 0}}{\mathbf{ \pi }}^{\mathbf{ 2}
\mathbf{ s}}_{\mathbf{ 0}}\mathbf{ -}{\mathbf{ \pi }}^{\mathbf{ 1}
\mathbf{ s}}_{\mathbf{ 0}}\sum_{i\in H}{\sum_{j\in H}{\sum_{\mathbf{ t}\in T}{p^s{\pi }^{st}_{2ij}X^{st}_{ij}}}}\mathbf{ -}{\mathbf{ \pi }}^{\mathbf{ 2}\mathbf{ s}}_{\mathbf{ 0}}\sum_{i\in H}{\sum_{t\in T}{p^s{\pi }^{st}_{1i}\left(Z^t_i-V^{t,t-1}_i\right)}}}}{{\sigma }_{{\pi }^{st}_{1}}{\sigma }_{{\pi }^{st}_{2}}}\right.\right.\nonumber\\
&\qquad\qquad\qquad\qquad\qquad\qquad  \left.\left.+\frac{\sum_{i\in H}{\sum_{j\in H}{\sum_{\mathbf{ t}\in T}{p^s{\pi }^{st}_{2ij}{\pi }^{st}_{1i}(L^{st}_{ij}-Q^{st,t-1}_{ij})}}}}{{\sigma }_{{\pi }^{st}_{1}}{\sigma }_{{\pi }^{st}_{2}}}-{\Psi }^*_s\right)-\gamma \right) \nonumber \\
&\qquad\qquad\qquad\qquad\qquad\qquad  +\sum_{s\in S}\sum_{t\in T}d_{2st}\left(1-\sum_{i\in H}\sum_{j\in H}X^{st}_{ij} \right)\nonumber
\end{align}
subject to: \eqref{eq5}, \eqref{eq16}-\eqref{eq21}
\[ 
X^{st}_{ij},\gamma \geq 0, Z^t_i, V^{t,t-1}_i, L^{st}_{ij}, Q^{st,t-1}_{ij}\in {\mathbb{ B}}^{\left|H\right|}\qquad i,j\in H, s\in S, t\in T\]
Where $d_1$and $d_2$ are vectors of Lagrange multipliers $d_{1s}$ and $d_{2st}$ for each $s\in S$ and $t\in T$, respectively.

\subsubsection{Benders decomposition algorithm}

Despite the removing the complicated constraints by LR method, the relaxed mixed integer program (RMIP) has a block structure in constraints \eqref{eq5}, \eqref{eq16}-\eqref{eq21}. Where, existence of $Z^t_i, V^{t,t-1}_i, L^{st}_{ij}, Q^{st,t-1}_{ij}\in {\mathbb{ B}}^{\left|H\right|}$ as complicating variables led to a dependency on blocks and this prevents separable solvability of the model. A decomposition approach such as BD method with single cut (SBD) uses this feature and can be capable to solve the large instances efficiently. The relaxed subproblem (SLRP) for given vectors ${\hat{Z}}^t_i, {\hat{V}}^{t,t-1}_i, {\hat{L}}^{st}_{ij}, {\hat{Q}}^{st,t-1}_{ij}\in {\mathbb{ B}}^{\left|H\right|}$ is written as follows:

$SLRP[d_{1} ,d_{2}]:$

\begin{align}
{\rm Min} \varphi= &{\theta }_1\left(\frac{\gamma }{{\gamma }^{\max}-{\gamma }^*}\right)+{\theta }_2\left(\frac{\sum_{t\in T}\sum_{i\in H}\sum_{j\in H}\sum_{s\in S}p^sw^{st}_i c^{st}_{ij}X^{st}_{ij}}{{\Omega }^{\max}-{\Omega }^*}\right)\label{eq32}\\
&\qquad\qquad +\sum_{s\in S}{d_{1s}}\left(\left(\frac{\mathbf{ -}{\mathbf{ \pi }}^{\mathbf{ 1}\mathbf{ s}}_{\mathbf{ 0}}\sum_{i\in H}{\sum_{j\in H}{\sum_{\mathbf{ t}\in T}{p^s{\pi }^{st}_{2ij}X^{st}_{ij}}}}}{{\sigma }_{{\pi }^{st}_{1}}{\sigma }_{{\pi }^{st}_{2}}}\right)-\gamma \right)\nonumber\\
&\qquad\qquad -\sum_{s\in S}{\sum_{t\in T}{d_{2st}\left(\sum_{i\in H}{\sum_{j\in H}{X^{st}_{ij}}} \right)}}\nonumber
\end{align}
Subject to:
\begin{align}
&\sum_{j\in H}{X^{st}_{ij}}+\sum_{j\in H/\{i\}}{X^{st}_{ji}}\leq {\hat{Z}}^t_i&&i\in H, s\in S, t\in T\label{eq33}\\
&X^{st}_{ij}\leq {\hat{L}}^{st}_{ij}-{\hat{Z}}^t_i+1&&i,j\in H, s\in S, t\in T\label{eq34}\\
&{-X}^{st}_{ij}\leq {\hat{Z}}^t_i-2{\hat{L}}^{st}_{ij}&&i,j\in H, s\in S, t\in T\label{eq35}\\
& X^{st}_{ij}\leq {\hat{Q}}^{st,t-1}_{ij}-{\hat{V}}^{t,t-1}_i+1&&i,j\in H, s\in S, t\in T\label{eq36}\\
&-X^{st}_{ij}\le {\hat{V}}^{t,t-1}_i-2{\hat{Q}}^{st,t-1}_{ij}&&i,j\in H, s\in S, t\in T\label{eq37}\\
&X^{st}_{ij}, \gamma \geq 0,&&i,j\in H, s\in S, t\in T\nonumber
\end{align}
The variables $u_{r} \geq 0$, $r=1,\ldots,5$ are defined as the dual variables associated to constraints \eqref{eq33}-\eqref{eq37}, respectively. Now the relaxed dual subproblem (DLRP) can be formulated to get upper bound (${\rm UB_{BD}}$) as follows:

$DLRP[d_{1} ,d_{2}]$:

\begin{align}
&{\rm Max}\ U=-\sum_{i\in H}{\sum_{s\in S}{\sum_{t\in T}{u^{st}_{1i}{\hat{Z}}^t_i}}}\mathbf{ -}\sum_{i\in H}{\sum_{j\in H}{\sum_{\mathbf{ t}\in T}{\sum_{s\in S}{u^{st}_{2ij}\left({\hat{L}}^{st}_{ij}-{\hat{Z}}^t_i+1\right)\ }}}}\label{eq38}\\
&-\sum_{i\in H}{\sum_{j\in H}{\sum_{\mathbf{t}\in T}{\sum_{s\in S}{u^{st}_{3ij}\left({\hat{Z}}^t_i-2{\hat{L}}^{st}_{ij}\right)\ }}}}-\sum_{i\in H}{\sum_{j\in H}{\sum_{\mathbf{ t}\in T}{\sum_{s\in S}{u^{st}_{4ij}\left({\hat{Q}}^{st,t-1}_{ij}-{\hat{V}}^{t,t-1}_i+1\right)}}}}\nonumber\\
&-\sum_{i\in H}{\sum_{j\in H}{\sum_{\mathbf{ t}\in T}{\sum_{s\in S}{u^{st}_{5ij}\left({\hat{V}}^{t,t-1}_i-2{\hat{Q}}^{st,t-1}_{ij}\right)\ }}}}\nonumber
\end{align}
Subject to:
\begin{align}
&-u^{st}_{1i}-u^{st}_{1j}-u^{st}_{2ij}+u^{st}_{3ij}-u^{st}_{4ij}+u^{st}_{5ij}\leq &\nonumber\\
&{\theta }_2\left(\frac{p^sw^{st}_i{\ c}^{st}_{ij}}{{\Omega }^{max}-{\Omega }^*}\right)-\frac{d_{1s}{\mathbf{ \pi }}^{\mathbf{ 1}\mathbf{ s}}_{\mathbf{ 0}}p^s\pi^{st}_{2ij}}{{\sigma }_{{\pi }^{st}_{1}}{\sigma }_{{\pi }^{st}_{2}}}-d_{2st}&i,j\in H, s\in S, t\in T\label{eq39}\\
&u_r\geq 0,\qquad r=1,\ldots, 5 &\label{eq40}
\end{align}
With the help of an auxiliary variable $\eta$ (as an under-estimator variable) in the corresponding objective function for the dual of optimality check subproblem, the Benders relaxed master problem (RMP) can be assembled as: $RMP [d_1,d_2]$:

\begin{align}\label{eq41}
{\rm Min}~\vartheta=&\eta-\frac{{{\theta }_1\gamma }^*}{{\gamma }^{\max}-{\gamma }^*}+{\theta }_2\left(\frac{\sum_{t\in T}{\sum_{i\in H}{f^t_i\left(Z^t_i-V^{t,t-1}_i\right)}}-{\Omega }^*}{{\Omega }^{max}-{\Omega }^*}\right)\\
&+\sum_{s\in S}{d_{1s}}\left(\frac{{\mathbf{ \pi }}^{\mathbf{ 1}\mathbf{ s}}_{\mathbf{ 0}}{\mathbf{ \pi }}^{\mathbf{ 2}\mathbf{ s}}_{\mathbf{ 0}}\mathbf{ -}{\mathbf{ \pi }}^{\mathbf{ 2}\mathbf{ s}}_{\mathbf{ 0}}\sum_{i\in H}{\sum_{t\in T}{p^s{\pi }^{st}_{1i}\left(Z^t_i-V^{t,t-1}_i\right)}}}{{\sigma }_{{\pi }^{st}_{1}}{\sigma }_{{\pi }^{st}_{2}}}\nonumber\right.\\
&\left.+\frac{\sum_{i\in H}{\sum_{j\in H}{\sum_{\mathbf{ t}\in T}{p^s{\pi }^{st}_{2ij}{\pi }^{st}_{1i}\left(L^{st}_{ij}-Q^{st,t-1}_{ij}\right)}}}}{{\sigma }_{{\pi }^{st}_{1}}{\sigma }_{{\pi }^{st}_{2}}}-{\Psi }^*_s\right)+\sum_{s\in S}{\sum_{t\in T}{d_{2st}}}\nonumber
\end{align}
Subject to: \eqref{eq16}-\eqref{eq17}
\[
Z^t_i, V^{t,t-1}_i, L^{st}_{ij}, Q^{st,t-1}_{ij}\in {\mathbb{ B}}^{\left|H\right|}\qquad i,j\in H, s\in S, t\in T\]

Let $\overline{u}_r\in E(p_D)$ and ${\overline{\overline{u}}}_r\in R\left(p_D\right)$, $r =1,\ldots ,5$ be the set of extreme points and extreme rays of the polyhedron of DLRP which is denoted by $p_D$, respectively. Now the relaxation of RMP (RRMP) can be formulated as follows:

$RRMP\left[d_1,d_2\right]$:

\begin{align*}
{\rm Min}~\vartheta & \\
=&\eta-\frac{{{\theta }_1\gamma }^*}{{\gamma }^{max}-{\gamma }^*}+{\theta }_2\left(\frac{\sum_{t\in T}{\sum_{i\in H}{f^t_i\left(Z^t_i-V^{t,t-1}_i\right)}}-{\Omega }^*}{{\Omega }^{max}-{\Omega }^*}\right)\\
&+\sum_{s\in S}{d_{1s}}\left(\frac{{\mathbf{ \pi }}^{\mathbf{ 1}\mathbf{ s}}_{\mathbf{ 0}}{\mathbf{ \pi }}^{\mathbf{ 2}\mathbf{ s}}_{\mathbf{ 0}}\mathbf{ -}{\mathbf{ \pi }}^{\mathbf{ 2}\mathbf{ s}}_{\mathbf{ 0}}\sum_{i\in H}{\sum_{t\in T}{\sum_{s\in S}{p^s{\pi }^{st}_{1i}\left(Z^t_i-V^{t,t-1}_i\right)}}}}{{\sigma }_{{\pi }^{st}_{1}}{\sigma }_{{\pi }^{st}_{2}}}\right.\\
&\left.+\frac{\sum_{i\in H}{\sum_{j\in H}{\sum_{\mathbf{ t}\in T}{\sum_{s\in S}{p^s{\pi }^{st}_{2ij}{\pi }^{st}_{1i}\left(L^{st}_{ij}-Q^{st,t-1}_{ij}\right)}}}}}{{\sigma }_{{\pi }^{st}_{1}}{\sigma }_{{\pi }^{st}_{2}}}-{\Psi }^*_s\right)+\sum_{s\in S}{\sum_{t\in T}{d_{2st}}}
\end{align*}
Subject to: \eqref{eq16}-\eqref{eq17}
\begin{align}\label{eq42}
\eta &\geq-\sum_{i\in H}{\sum_{t\in T}{\sum_{s\in S}{{\overline{u}}^{st}_{1i}Z^t_i}}}-\sum_{i\in H}{\sum_{j\in H}{\sum_{t\in T}{\sum_{s\in S}{{\overline{u}}^{st}_{2ij}\left(L^{st}_{ij}-Z^t_i+1\right)\ }}}}\nonumber\\
&-\sum_{i\in H}{\sum_{j\in H}{\sum_{t\in T}{\sum_{s\in S}{{\overline{u}}^{st}_{3ij}\left(Z^t_i-2L^{st}_{ij}\right)\ }}}}
-\sum_{i\in H}{\sum_{t\in T}{\sum_{\mathbf{ j}\in H}{\sum_{s\in S}{{\overline{u}}^{st}_{4ij}\left(Q^{st,t-1}_{ij}-V^{t,t-1}_i+1\right)\ }}}}\nonumber\\
&-\sum_{i\in H}{\sum_{t\in T}{\sum_{\mathbf{ j}\in H}{\sum_{s\in S}{{\overline{u}}^{st}_{5ij}\left(V^{t,t-1}_i-2Q^{st,t-1}_{ij}\right)\ }}}}\qquad \forall {\overline{u}}_r\in E\left(p_D\right)
\end{align}

\begin{align}
0&\geq -\sum_{i\in H}\sum_{t\in T}\sum_{s\in S}\overline{\overline{u}}^{st}_{1i}Z^t_i-\sum_{i\in H}{\sum_{j\in H}{\sum_{t\in T}{\sum_{s\in S}{{\overline{\overline{u}}}^{st}_{2ij}\left(L^{st}_{ij}-Z^t_i+1\right)}}}}\nonumber\\
&-\sum_{i\in H}\sum_{j\in H}\sum_{t\in T}\sum_{s\in S}{\overline{\overline{u}}}^{st}_{3ij}\left(Z^t_i-2L^{st}_{ij}\right)
-\sum_{i\in H}{\sum_{t\in T}{\sum_{\mathbf{ j}\in H}{\sum_{s\in S}{{\overline{\overline{u}}}^{st}_{4ij}\left(Q^{st,t-1}_{ij}-V^{t,t-1}_i+1\right)}}}}\nonumber\\
&-\sum_{i\in H}{\sum_{t\in T}{\sum_{\mathbf{ j}\in H}{\sum_{s\in S}{{\overline{\overline{u}}}^{st}_{5ij}\left(V^{t,t-1}_i-2Q^{st,t-1}_{ij}\right)}}}}\qquad\forall {\overline{\overline{u}}}_r\in R(p_D)
\end{align} 
\[Z^t_i,V^{t,t-1}_i,L^{st}_{ij},Q^{st,t-1}_{ij}\in {\mathbb{ B}}^{\left|H\right|}\qquad i,j\in H, s\in S, t\in T\]

\subsubsection{BD enhancement}

Because of the structure of LRP, during the decomposition of the problem many of constraints (i.e. Eqs, \eqref{eq5}, \eqref{eq16}-\eqref{eq21}) appeared in the SP with a resulting lack of effect on the master problem. Therefore, adding only one cut at each iteration will result slow convergence and a large gap associated with the classical BD algorithm. Benders Decomposition with multi-cuts should be an efficient alternative. Moreover, the network structure of the problem has the consequence of high degeneracy being exhibited in the SP and there is also the possibility of alternative cuts. In this situation, generating strong Pareto cuts plays an important role in enhancing the algorithms efficiency. In this way, we are motivated to employ the multi-Pareto cut BD algorithm (MPBD) to solve LRP (LRMPBD) with less computational effort.

\paragraph{Multiple cut generation strategies}

The independency of the solution space of DSP on complicating variables allows us to generate all cuts in the first iteration whereas the more complexity in RMP constraints is not desirable. Hence, there is a tradeoff between the number of cuts and iterations that should be considered. Here, Multi-cut generation (called MBD) can play a crucial role as a strategy to accelerate convergence. In our case SP can be decomposed into $|S|\times |T|$ independent SPs  and we generate one cut only for each subproblem in each iteration. We may reformulate RRMP (MRRMP) by defining $\eta^{st}$ and $p_{D^{st}}$ to be the network total cost and the dual polyhedral of each SPs of period time $t\in T$, based on scenario $s\in S$, as follows.

\begin{align}\label{eq44}
{\rm Min}\ \vartheta=&\sum_{s\in S}{\sum_{t\in T}{{\rm \ }\eta^{st}}}-\frac{{{\theta }_1\gamma }^*}{{\gamma }^{max}-{\gamma }^*}+{\theta }_2\left(\frac{\sum_{t\in T}{\sum_{i\in H}{f^t_i\left(Z^t_i-V^{t,t-1}_i\right)}}-{\Omega }^*}{{\Omega }^{max}-{\Omega }^*}\right)\\
&+\sum_{s\in S}{d_{1s}}\left(\frac{{\mathbf{ \pi }}^{\mathbf{ 1}\mathbf{ s}}_{\mathbf{ 0}}{\mathbf{ \pi }}^{\mathbf{ 2}\mathbf{ s}}_{\mathbf{ 0}}\mathbf{ -}{\mathbf{ \pi }}^{\mathbf{ 2}\mathbf{ s}}_{\mathbf{ 0}}\sum_{i\in H}{\sum_{t\in T}{\sum_{s\in S}{p^s{\pi }^{st}_{1i}\left(Z^t_i-V^{t,t-1}_i\right)}}}}{{\sigma }_{{\pi }^{st}_{1}}{\sigma }_{{\pi }^{st}_{2}}}\right.\nonumber\\
&\left.+\frac{\sum_{i\in H}{\sum_{j\in H}{\sum_{\mathbf{ t}\in T}{\sum_{s\in S}{p^s{\pi }^{st}_{2ij}{\pi }^{st}_{1i}\left(L^{st}_{ij}-Q^{st,t-1}_{ij}\right)}}}}}{{\sigma }_{{\pi }^{st}_{1}}{\sigma }_{{\pi }^{st}_{2}}}-{\Psi }^*_s\right)+\sum_{s\in S}{\sum_{t\in T}{d_{2st}}}\nonumber
\end{align}
Subject to: \eqref{eq16}-\eqref{eq17}
\begin{align}\label{eq45}
\eta^{st}\geq &-\sum_{i\in H}{{\overline{u}}^{st}_{1i}Z^t_i}-\sum_{i\in H}{\sum_{j\in H}{{\overline{u}}^{st}_{2ij}\left(L^{st}_{ij}-Z^t_i+1\right)}}-\sum_{i\in H}{\sum_{j\in H}{{\overline{u}}^{st}_{3ij}\left(Z^t_i-2L^{st}_{ij}\right)}}\\
&-\sum_{i\in H}{\sum_{\mathbf{ j}\in H}{{\overline{u}}^{st}_{4ij}\left(Q^{st,t-1}_{ij}-V^{t,t-1}_i+1\right)}}\nonumber\\
&-\sum_{i\in H}{\sum_{\mathbf{ j}\in H}{{\overline{u}}^{st}_{5ij}\left(V^{t,t-1}_i-2Q^{st,t-1}_{ij}\right)}}\qquad \forall {\overline{u}}_r\in E\left(p_{D^{st}}\right),\ s\in S\ \ ,\ t\in T\nonumber
\end{align}

\begin{align}\label{eq46}
0&\geq -\sum_{i\in H}{{\overline{\overline{u}}}^{st}_{1i}Z^t_i}\mathbf{ -}\sum_{i\in H}{\sum_{j\in H}{{\overline{\overline{u}}}^{st}_{2ij}\left(L^{st}_{ij}-Z^t_i+1\right)}}\mathbf{ -}\sum_{i\in H}{\sum_{j\in H}{{\overline{\overline{u}}}^{st}_{3ij}\left(Z^t_i-2L^{st}_{ij}\right)}}\\
&-\sum_{i\in H}{\sum_{\mathbf{ j}\in H}{{\overline{\overline{u}}}^{st}_{4ij}\left(Q^{st,t-1}_{ij}-V^{t,t-1}_i+1\right)}}\mathbf{ -}\sum_{i\in H}{\sum_{\mathbf{ j}\in H}{{\overline{\overline{u}}}^{st}_{5ij}\left(V^{t,t-1}_i-2Q^{st,t-1}_{ij}\right)}}\nonumber\\
&\hspace*{6.8cm}\forall {\overline{\overline{u}}}_r\in R(p_{D^{st}}),\ s\in S\ \ ,\ t\in T\nonumber
\end{align}

\[Z^t_i, V^{t,t-1}_i, L^{st}_{ij}, Q^{st,t-1}_{ij}\in {\mathbb{ B}}^{\left|H\right|}\qquad i,j\in H, s\in S, t\in T\]

Noteworthy, solving the MRRMP problem achieve the ${\rm LB_{BD}}$ which is considered as the ${\rm LB_{LR}}$ along the solution algorithm.

\paragraph{Pareto cut generation strategies}

Because of the network structure of HLPs, PRH-R, typically we have multiple dual solutions which will have consequences for the set of optimal cuts. Here, the Pareto cut strategy (called PBD) can be helpful to construct a stronger cut. The generated cut corresponding to a solution ${\overline{u}}^a_r$ dominates a generated cut which is associated to solution ${\overline{u}}^b_r$ if and only if:

\begin{align}\label{eq47}
&-\sum_{i\in H}\sum_{s\in S}\sum_{t\in T}{\overline{u}^{st}_{1i}}^aZ^t_i-\sum_{i\in H}\sum_{j\in H}\sum_{t\in T}\sum_{s\in S}{\overline{u}^{st}_{2ij}}^a\left(L^{st}_{ij}-Z^t_i+1\right)\cr
&-\sum_{i\in H}\sum_{j\in H}\sum_{t\in T}\sum_{s\in S}{\overline{u}^{st}_{3ij}}^a\left(Z^t_i-2L^{st}_{ij}\right)\cr
&-\sum_{i\in H}{\sum_{j\in H}{\sum_{\mathbf{ t}\in T}{\sum_{s\in S}{{{\overline{u}}^{st}_{4ij}}^a\left(Q^{st,t-1}_{ij}-V^{t,t-1}_i+1\right)\ }}}}\\
&\mathbf{ -}\sum_{i\in H}{\sum_{j\in H}{\sum_{\mathbf{ t}\in T}{\sum_{s\in S}{{{\overline{u}}^{st}_{5ij}}^a\left(V^{t,t-1}_i-2Q^{st,t-1}_{ij}\right)\ }}}}\nonumber\\
&\geq -\sum_{i\in H}{\sum_{s\in S}{\sum_{t\in T}{{{\overline{u}}^{st}_{1i}}^bZ^t_i}}}-\sum_{i\in H}{\sum_{j\in H}{\sum_{\mathbf{ t}\in T}{\sum_{s\in S}{{{\overline{u}}^{st}_{2ij}}^b\left(L^{st}_{ij}-Z^t_i+1\right)\ }}}}\nonumber\\
&-\sum_{i\in H}{\sum_{j\in H}{\sum_{\mathbf{ t}\in T}{\sum_{s\in S}{{{\overline{u}}^{st}_{3ij}}^b\left(Z^t_i-2L^{st}_{ij}\right)\ }}}}-\sum_{i\in H}{\sum_{j\in H}{\sum_{\mathbf{ t}\in T}{\sum_{s\in S}{{{\overline{u}}^{st}_{4ij}}^b\left(Q^{st,t-1}_{ij}-V^{t,t-1}_i+1\right)\ }}}}\nonumber\\
&-\sum_{i\in H}{\sum_{j\in H}{\sum_{\mathbf{ t}\in T}{\sum_{s\in S}{{{\overline{u}}^{st}_{5ij}}^b\left(V^{t,t-1}_i-2Q^{st,t-1}_{ij}\right)\ }}}}\qquad\forall {\overline{u}}_r\in E\left(p_D\right)\nonumber
\end{align}

The Pareto cut is a cut which dominates all other cuts. The Pareto cut sets are generated by solving the following auxiliary dual problem:

\begin{align}\label{eq48}
{\rm Max}\ Prt=&-\sum_{i\in H}{\sum_{s\in S}{\sum_{t\in T}{u^{st}_{1i}{\dot{Z}}^t_i}}}\mathbf{ -}\sum_{i\in H}{\sum_{j\in H}{\sum_{\mathbf{ t}\in T}{\sum_{s\in S}{u^{st}_{2ij}\left({\dot{L}}^{st}_{ij}-{\dot{Z}}^t_i+1\right)\ }}}}\\
&-\sum_{i\in H}{\sum_{j\in H}{\sum_{\mathbf{ t}\in T}{\sum_{s\in S}{u^{st}_{3ij}\left({\dot{Z}}^t_i-2{\dot{L}}^{st}_{ij}\right)\ }}}}\cr
&-\sum_{i\in H}{\sum_{j\in H}{\sum_{\mathbf{ t}\in T}{\sum_{s\in S}{u^{st}_{4ij}\left({\dot{Q}}^{st,t-1}_{ij}-{\dot{V}}^{t,t-1}_i+1\right)\ }}}}\nonumber\\
&-\sum_{i\in H}{\sum_{j\in H}{\sum_{\mathbf{ t}\in T}{\sum_{s\in S}{u^{st}_{5ij}\left({\dot{V}}^{t,t-1}_i-2{\dot{Q}}^{st,t-1}_{ij}\right)\ }}}}\nonumber
\end{align} 

Subject to: \eqref{eq39}
\begin{align}\label{eq49}
&-\sum_{i\in H}{\sum_{s\in S}{\sum_{t\in T}{u^{st}_{1i}{\hat{Z}}^t_i}}}\mathbf{ -}\sum_{i\in H}{\sum_{j\in H}{\sum_{\mathbf{ t}\in T}{\sum_{s\in S}{u^{st}_{2ij}\left({\hat{L}}^{st}_{ij}-{\hat{Z}}^t_i+1\right)\ }}}}\cr
&\mathbf{ -}\sum_{i\in H}{\sum_{j\in H}{\sum_{\mathbf{ t}\in T}{\sum_{s\in S}{u^{st}_{3ij}\left({\hat{Z}}^t_i-2{\hat{L}}^{st}_{ij}\right)\ }}}}\\
&\mathbf{ -}\sum_{i\in H}{\sum_{j\in H}{\sum_{\mathbf{ t}\in T}{\sum_{s\in S}{u^{st}_{4ij}\left({\hat{Q}}^{st,t-1}_{ij}-{\hat{V}}^{t,t-1}_i+1\right)\ }}}}\cr
&\mathbf{ -}\sum_{i\in H}{\sum_{j\in H}{\sum_{\mathbf{ t}\in T}{\sum_{s\in S}{u^{st}_{5ij}\left({\hat{V}}^{t,t-1}_i-2{\hat{Q}}^{st,t-1}_{ij}\right)\ }}}}\mathbf{ =}{{ U}}^{\mathbf{ *}}\nonumber\\
&\hspace*{8.5cm} u_r\geq 0,\quad r=1,\ldots,5\nonumber
\end{align}

where, $U^*$ is the optimal solution of DLRP. Also, $\dot{Z},\dot{V},\dot{L}$ and $\dot{Q}$ refer to core points i.e. a set of interior points of RMP convex hull. 

Papadakos (2008) showed that the convergence rate is affected by the value of core points, so a linear combination of current core points and their values corresponding to the latest iteration is considered as an intensification procedure.

\paragraph{Multi-Pareto cut generation strategies}

Generating the large number of cuts led to complexity in RMP. Therefore, a multi-Pareto strategy can improve the efficiency of the algorithm. We reformulate the auxiliary dual problem (ADP) as follows:

\begin{align}\label{eq50}
\mathop{\max Prt^{st}}_{s\in S,\ t\in T}&=-\sum_{i\in H}{u^{st}_{1i}{\dot{Z}}^t_i}-\sum_{i\in H}{\sum_{j\in H}{u^{st}_{2ij}\left({\dot{L}}^{st}_{ij}-{\dot{Z}}^t_i+1\right)}}\\
&-\sum_{i\in H}{\sum_{j\in H}{u^{st}_{3ij}\left({\dot{Z}}^t_i-2{\dot{L}}^{st}_{ij}\right)}}-\sum_{i\in H}{\sum_{j\in H}{u^{st}_{4ij}\left({\dot{Q}}^{st,t-1}_{ij}-{\dot{V}}^{t,t-1}_i+1\right)}}\nonumber\\
&-\sum_{i\in H}{\sum_{j\in H}{u^{st}_{5ij}\left({\dot{V}}^{t,t-1}_i-2{\dot{Q}}^{st,t-1}_{ij}\right)}}\nonumber
\end{align}

Subject to: \eqref{eq39}

\begin{align}\label{eq51}
&-\sum_{i\in H}{u^{st}_{1i}{\hat{Z}}^t_i}-\sum_{i\in H}{\sum_{j\in H}{u^{st}_{2ij}\left({\hat{L}}^{st}_{ij}-{\hat{Z}}^t_i+1\right)}}-\sum_{i\in H}{\sum_{j\in H}{u^{st}_{3ij}\left({\hat{Z}}^t_i-2{\hat{L}}^{st}_{ij}\right)}}\\
&-\sum_{i\in H}{\sum_{j\in H}{u^{st}_{4ij}\left({\hat{Q}}^{st,t-1}_{ij}-{\hat{V}}^{t,t-1}_i+1\right)}}\mathbf{ -}\sum_{i\in H}{\sum_{j\in H}{u^{st}_{5ij}\left({\hat{V}}^{t,t-1}_i-2{\hat{Q}}^{st,t-1}_{ij}\right)}}\nonumber\\
&={U^*}^{st}\qquad s\in S, t\in T\nonumber\\
&u_r \geq 0,\quad r=1,\ldots ,5\nonumber
\end{align}

Where, the optimal solution of DLRP for scenario $s\in S$ in period time $t\in T$ is symbolized by ${U^*}^{st}$.

\subsubsection{Subgradient multipliers updating procedure}

The following subgradient optimization procedure is used to update the Lagrange multipliers iteratively (Fisher 2004). Denote by $d^{{it}_{LR}}_{1s}$ and $d^{{it}_{LR}}_{2st}$ the updated Lagrange multipliers vectors in each iteration$\ it_{LR}$. Also, by ${\theta }^{{it}_{LR}}$ and ${\sigma }^{{it}_{LR}}$ denote the ${it}_{LR}$-\textit{th} iteration step size and step size parameter, respectively. Then the step size can be calculated as follows:

\begin{equation}\label{eq52}
{\theta }^{{it}_{LR}}={\sigma }^{{it}_{LR}}\frac{{UB}^{{it}_{LR}}_{LR}-{LB}_{BD}{\rm (}d^{it_{LR}}_1,d^{{it}_{LR}}_2{\rm )}}{{\|{\triangle d}^{it_{LR}}_1,{\triangle d}^{it_{LR}}_2\|}^2}
\end{equation}

Here ${LB}_{BD}(d^{{it}_{LR}}_{1s},d^{{it}_{LR}}_{2st})$ provides the BD lower bound corresponding to given Lagrange multipliers $d^{{it}_{LR}}_{1s}$ and $d^{{it}_{LR}}_{2st}$ in iteration ${it}_{LR}$. Also, the subgradient vector associated with constraints \eqref{eq5} and \eqref{eq14} can be determined as:

{\small
\begin{align}\label{eq53}
&{\left|\left|{\triangle d}^{it_{LR}}_1,{\triangle d}^{it_{LR}}_2\right|\right|}^2=\\ 
&\sum_{s\in S}\left(\left(\frac{{\mathbf{ \pi }}^{\mathbf{ 1}\mathbf{ s}}_{\mathbf{ 0}}{\mathbf{ \pi }}^{\mathbf{ 2}\mathbf{ s}}_{\mathbf{ 0}}\mathbf{ -}{\mathbf{ \pi }}^{\mathbf{ 1}\mathbf{ s}}_{\mathbf{ 0}}\sum_{i\in H}{\sum_{j\in H}{\sum_{\mathbf{ t}\in T}{\sum_{s\in S}{p^s{\pi }^{st}_{2ij}{\hat{X}}^{st}_{ij}}}}}\mathbf{ -}{\mathbf{ \pi }}^{\mathbf{ 2}\mathbf{ s}}_{\mathbf{ 0}}\sum_{i\in H}{\sum_{t\in T}{\sum_{s\in S}{p^s{\pi }^{st}_{1i}\left({\hat{Z}}^t_i-{\hat{V}}^{t,t-1}_i\right)}}}}{{\sigma }_{{\pi }^{st}_{1}}{\sigma }_{{\pi }^{st}_{2}}}\right.\right.\nonumber\\
&\qquad\qquad\qquad  +\left.\left.\frac{\sum_{i\in H}{\sum_{j\in H}{\sum_{\mathbf{ t}\in T}{\sum_{s\in S}{p^s{\pi }^{st}_{2ij}{\pi }^{st}_{1i}\left({\hat{L}}^{st}_{ij}-{\hat{Q}}^{st,t-1}_{ij}\right)}}}}}{{\sigma }_{{\pi }^{st}_{1}}{\sigma }_{{\pi }^{st}_{2}}}-{\Psi }^*_s\right)-\widehat{\gamma }\right)^2 \nonumber \\
&\qquad\qquad\qquad +\sum_{t\in T}\left(1-\sum_{i\in H}{\sum_{j\in H}{{\hat{X}}^{st}_{ij}}}\right)^2\nonumber 
\end{align}
}

It is worth to note that UB${}_{LR}$ denotes an upper bound for the original problem \eqref{eq23}, \eqref{eq4}-\eqref{eq5}, \eqref{eq15}-\eqref{eq22}. To obtain this, the original problem is solved with fixed values of the variable set $V^{t,t-1}_i$. Fixing value of this variable produces a solvable problem which result in an efficient upper bound whereas fixing value of other variables lead to a hard restricted model.

The Lagrange multipliers are updated as follows:

\begin{align}\label{eq54} 
&d^{t_{LR}}_{1s}=\\
&\max\left\{0,d^{t_{{LR}^{-1}}}_{1s}\right.\nonumber\\
&+{\theta }^{t_{{LR}^{-1}}}\left(\left(\frac{{\mathbf{ \pi }}^{\mathbf{ 1}\mathbf{ s}}_{\mathbf{ 0}}{\mathbf{ \pi }}^{\mathbf{ 2}\mathbf{ s}}_{\mathbf{ 0}}\mathbf{ -}{\mathbf{ \pi }}^{\mathbf{ 1}\mathbf{ s}}_{\mathbf{ 0}}\sum_{i\in H}{\sum_{j\in H}{\sum_{\mathbf{ t}\in T}{\sum_{s\in S}{p^s{\pi }^{st}_{2ij}{\hat{X}}^{st}_{ij}}}}}\mathbf{ -}{\mathbf{ \pi }}^{\mathbf{ 2}\mathbf{ s}}_{\mathbf{ 0}}\sum_{i\in H}{\sum_{t\in T}{\sum_{s\in S}{p^s{\pi }^{st}_{1i}\left({\hat{Z}}^t_i-{\hat{V}}^{t,t-1}_i\right)}}}}{{\sigma }_{{\pi }^{st}_{1}}{\sigma }_{{\pi }^{st}_{2}}}\right.\right.\nonumber\\
&+\left.\left.\left.\frac{\sum_{i\in H}{\sum_{j\in H}{\sum_{\mathbf{ t}\in T}{\sum_{s\in S}{p^s{\pi }^{st}_{2ij}{\pi }^{st}_{1i}({\hat{L}}^{st}_{ij}-{\hat{Q}}^{st,t-1}_{ij})}}}}}{{\sigma }_{{\pi }^{st}_{1}}{\sigma }_{{\pi }^{st}_{2}}}-{\Psi }^*_s\right)-\widehat{\gamma}\right)\right\}\qquad s\in S\nonumber
\end{align}
\begin{equation}\label{eq55}
d^{t_{LR}}_{2st}=\max\left\{0,d^{t_{{LR}^{-1}}}_{2st}+{\theta }^{t_{{LR}^{-1}}}\left(1-\sum_{i\in H}\sum_{j\in H}\hat{X}^{st}_{ij}\right)\right\}
\qquad s\in S, t\in T
\end{equation}

Updating the Lagrange multipliers result in reducing the gap and terminate the algorithm, iteratively.

Let ${it}_{BD}$ represents the inner iteration number of the BD method. Now, the pseudo-code of the proposed relax-and-decomposition algorithm is depicted in Algorithm 1.

\begin{scriptsize}
\begin{tabular}{p{4.9in}} \hline 
\textbf{Algorithm 1}: Relaxed-and-Decomposition heuristic algorithm for \textit{PRH-R} \\ \hline 
set $d^0_{1s}=d^0_{2st}=0$, ${\sigma }^0$\textit{=}2, \textit{noipr}=0, ${LB}^0_{LR}$=${LB}^0_{BD}$=-âˆž,$\ {UB}^0_{LR}={UB}^0_{BD}=\infty $\newline \textbf{repet} \\ \hline 
\textit{Solving LRP by the Multi-Pareto cut BD method} \\ \hline 

set \textit{initial core points} $\dot{Z},\dot{V},\dot{L},\dot{Q}$\newline 
${\rm solve\ }DLRP\left({\hat{Z}}^0,{\hat{V}}^0,{\hat{L}}^0,{\hat{Q}}^0|{\hat{d}}^{{it}_{LR}}_{1s},{\hat{d}}^{it_{LR}}_{2st}\right)$\textit{ for each scenario and period time}\newline 
\textbf{  if}\textit{  }$\ SLRP\left({\hat{Z}}^0,{\hat{V}}^0,{\hat{L}}^0,{\hat{Q}}^0|{\hat{d}}^{it_{LR}}_{1s},{\hat{d}}^{it_{LR}}_{2st}\right)\ infesible{\rm \ }then{\rm \ }obtain\ {\overline{\overline{u}}}^0_r$\textit{\newline }
\textbf{       }add\textit{ feasibility cut set (46) with} ${\overline{\overline{u}}}^0_r\ to\ MRRMP\ to\ obtain\ {\hat{Z}}^{{it}_{BD}},{\hat{V}}^{it_{BD}},{\hat{L}}^{it_{BD}},{\hat{Q}}^{it_{BD}},\ {it}_{BD}\ge 1$\newline
 \textbf{   else if} $SLRP\left({\hat{Z}}^0,{\hat{V}}^0,{\hat{L}}^0,{\hat{Q}}^0|{\hat{d}}^{it_{LR}}_{{\rm 1s}},{\hat{d}}^{it_{LR}}_{2st}\right) $ feasible\ and \newline 
$ \left({UB}_{BD} ({\hat{Z}}^0,{\hat{V}}^0,{\hat{L}}^0,{\hat{Q}}^0|{\hat{d}}^{it_{LR}}_{{\rm 1s}},{\hat{d}}^{it_{LR}}_{{\rm 2st}} ) -{LB}_{BD}{\rm (}{\hat{Z}}^0,{\hat{V}}^0,{\hat{L}}^0,{\hat{Q}}^0|{\hat{d}}^{it_{LR}}_{{\rm 1s}},{\hat{d}}^{it_{LR}}_{{\rm 2st}}{\rm )} \right)\le {\varepsilon }_{BD}\ \mathbf{ then}$\newline 
 
 \textbf{         stop\textit{\newline }       Else\newline }          update\textit{ core points }$\dot{Z},\dot{V},\dot{L},\dot{Q}$\newline 
         ${\rm solve}\ ADP\ to\ obtain\ {\overline{u}}^0_r$\textbf{\textit{\newline }}\textit{         }add\textit{ multi-Pareto optimality cut set (45) with }${\overline{u}}^0_r\ \ to\ MRRMP$\newline \textbf{         }solve $MRRMP({\overline{u}}^{{it}_{BD}}_r)\ \ to\ obtain\ {\hat{Z}}^{{it}_{BD}},{\hat{V}}^{it_{BD}},{\hat{L}}^{it_{BD}},{\hat{Q}}^{it_{BD}}$\newline \textbf{         }solve $DLRP\left({\hat{Z}}^{{it}_{BD}},{\hat{V}}^{it_{BD}},{\hat{L}}^{it_{BD}},{\hat{Q}}^{it_{BD}}|{\hat{d}}^{{it}_{LR}}_{1s},{\hat{d}}^{it_{LR}}_{2st}\right)\ $\textit{for each scenario and period time}\newline 
          \textbf{   if} $\ SLRP\left({\hat{Z}}^{{it}_{BD}},{\hat{V}}^{it_{BD}},{\hat{L}}^{it_{BD}},{\hat{Q}}^{it_{BD}}|{\hat{d}}^{{it}_{LR}}_{1s},{\hat{d}}^{it_{LR}}_{2st}\right)\ \ infesible\ then\ obtain\ {\overline{\overline{u}}}^{it_{BD}}_r$
           \newline       
        add \ \textit{ feasibility cut set (46) with }${\overline{\overline{u}}}^{it_{BD}}_r\ to\ MRRMP\ to\ obtain\ {\hat{Z}}^{{it}_{BD}},{\hat{V}}^{it_{BD}},{\hat{L}}^{it_{BD}},{\hat{Q}}^{it_{BD}}$\newline
         \textbf{  else if} $SLRP\left({\hat{Z}}^{{it}_{BD}},{\hat{V}}^{it_{BD}},{\hat{L}}^{it_{BD}},{\hat{Q}}^{it_{BD}}|{\hat{d}}^{{it}_{LR}}_{1s},{\hat{d}}^{it_{LR}}_{2st}\right)$ {\textit feasible and } \newline 
          $\left({UB}_{BD}\left({\hat{Z}}^{{it}_{BD}},{\hat{V}}^{it_{BD}},{\hat{L}}^{it_{BD}},{\hat{Q}}^{it_{BD}}|{\hat{d}}^{{it}_{LR}}_{1s},{\hat{d}}^{it_{LR}}_{2st}\right) \right.$ \newline 
       $\left. \qquad\qquad\qquad\qquad  {-}{LB}_{BD}\left({\hat{Z}}^{{it}_{BD}},{\hat{V}}^{it_{BD}},{\hat{L}}^{it_{BD}},{\hat{Q}}^{it_{BD}}|{\hat{d}}^{{it}_{LR}}_{1s},{\hat{d}}^{it_{LR}}_{2st}\right)\right)\le {\varepsilon }_{BD}\ \mathbf{ then}$\textbf{\newline            stop\textit{\newline }        Else
      	\newline }         
       update\textit{ core points }$\dot{Z},\dot{V},\dot{L},\dot{Q}$\textit{\newline }\textbf{        }${\rm solve}\ ADP\ to\ obtain\ {\overline{u}}^{{it}_{BD}}_r$\textbf{\textit{\newline }}\textit{         }add\textit{ multi-Pareto optimality cut set (45) with }${\overline{u}}^{{it}_{BD}}_r\ to\ MRRMP$\newline \textbf{     end\textit{\newline } end \textit{\newline }}\textit{   }${it}_{BD}$\textit{â†’ }${it}_{BD}$\textit{+1}\newline 
       
       \textbf{end }(if MPBD's stoping criteria are met)\textbf{} \\ \hline 
\textbf{    if} ${LB}_{BD}$($d^{{it}_{LR}}_{1s},d^{{it}_{LR}}_{2st}$) $>$ ${LB}^{it_{LR}}_{LR}$\newline \textbf{      then} set ${LB}^{it_{LR}}_{LR}$=${LB}_{BD}$ ($d^{it_{LR}}_{1s},d^{it_{LR}}_{2st}$), \textit{noipr}=0\newline \textbf{    else \newline }    increment \textit{noipr}\newline \textbf{        if} \textit{noipr} $>$ 1\newline           \textbf{then} set ${\sigma }^{it_{LR}}$= ${\sigma }^{{it}_{LR}-1}$/2, \textit{noipr}=0\newline         solve ${UB}^{it_{LR}}_{LR}$\textit{by obtained fixed variable} ${\overline{V}}^{it_{LR}}$ at ${LB}_{BD}$($d^{it_{LR}}_{1s},d^{it_{LR}}_{2st}$)\newline         calculate \textit{the subgradiant} $\triangle $\textit{:}\newline         calculate \textit{the step size} ${\theta }^{{it}_{LR}}={\sigma }^{{it}_{LR}}\frac{{UB}^{it_{LR}}_{LR}-{LB}^{it_{LR}}_{LR}}{{||{\triangle d}^{it_{LR}}_1,{\triangle d}^{{it}_{LR}}_2||}^2}$\newline         update \textit{multiplayeirs} $d^{{it}_{LR}}_{1s},d^{{it}_{LR}}_{2st}$\newline \textbf{until} (subgradient's stoping conditions are satisfied) \\ \hline 
\end{tabular}
\end{scriptsize}

The algorithm is repeated until stopping conditions are met, where the subgradient algorithm stopping conditions are considered as follows:

\begin{itemize}
\item  The step size parameter ${\sigma }^{{it}_{LR}}$ to be less than a prespecific threshold value ${\varepsilon }_{LR}$.  

\item  Reaching a maximum number of iterations $\left({Iter1}_{\max}\right)$.

\item  Reaching a maximum CPU time $\left({Time}_{\max}\right)$.
\end{itemize}

Also following stopping criteria have been considered for MPBD algorithm:

\begin{itemize}
\item  The percentage gap between the lower and upper bounds to be less than a threshold value ${\varepsilon }_{BD}$.

\item  Reaching a maximum number of iterations $\left({Iter2}_{max}\right)$.
\end{itemize}

\section{Experiment results}\label{se5}

We discussed on a series of computational experiments to evaluate the performance of the proposed model and the proposed relax-and-decomposition algorithm. First, the models validity is considered for a determined number of scenarios obtained by the SAA scheme. Then, the solution method behavior is analyzed. All numerical experiments are solved by an Asus Studio PC with an Intel Core i7 CPU at 1.73 GHz and 4 GB of RAM. Also, GAMS 23.5 optimization software and CPLEX solver are used to code the heuristic algorithm.

Numerical experiments are constructed based on the well-known CAB data set with minor modifications. The discount factor $\tau $ and the objective function weights are set to 0.8, 0.4 and 0.6, respectively. The stochastic flow between hub nodes is generated in different scenarios and period times according to a uniform distribution of [0.5, 0.8] $\times$ (distance) while distances are taken from the CAB data set. The stochastic unit transportation cost for different scenarios and period times are assumed to be uniformly distributed in [10, 20]. The setup cost of hub facilities is considered to be uniformly distributed in $[100\times \log  O_1 (i,t), 200\times\log  O_1\left(i,t\right)]$, where, $O_1\left(i,t\right)=U[0.5,1]\times\sum_j{w' \left(i, j\right)}$ and the flow of CAB data set is denoted by $w' \left(i\ ,\ j\right)$. The threshold value of risk for potential hub nodes and hub links are defined to be uniformly distributed in [\textbar \textit{H}\textbar \textbar \textit{T}\textbar ,10\textbar \textit{H}\textbar \textbar \textit{T}\textbar \textbar \textit{S}\textbar ]. After some tuning, values related to stopping criteria are set as: ${Iter1}_{max}=30,\ {Time}_{max}=10000(sec)$,$\ {Iter2}_{max}=20$, ${\varepsilon }_{LR}=0.001$ and ${\varepsilon }_{BD}=0.01$. Also, the core point update rule is defined according to Papadakos (2008). For example, it is calculated for $Z$ as $\dot{Z}=\lambda \dot{Z}+\left(1-\lambda \right)\hat{Z}$ where, $\lambda =0.5$ and core point sets related to $Z$ and $V$ are initialized to zero and core point sets corresponding to $L$ and $Q$ are initialized with a set of feasible solution to the RMP.

Some symbols which are used in the tables are as follows:

\begin{itemize}
\item
$|H|$: The number of total hub nodes.
\item
$|T|$: The number of period times.
\item
${LB}_{(.)}$: The lower bound for different solution approaches.
\item  ${\%Dev}_{DifLR}$: The percent deviation correspond each lower bound deal with the classical LR algorithm and the proposed solution method with respect to the best found lower bound.
\item  ${\%Dev}_{BD}$: The percent deviation correspond each lower bound deal with different version of BD algorithms with respect to the best found lower bound.
\item  $\%optimality~ gap$: The percent deviation between the upper and lower bound of the MPBD algorithm. That is $\%optimality~ gap=100(UB_{BD}$-$LB_{BD})/ UB_{BD}$.
\end{itemize}

\subsection{Practical convergence of SAA algorithm}

By increasing the sample size and the number of replications, the quality of the solution is increased as well as the computational complexity. For this, an efficient sample size is chooses making the trade-off between the accuracy of solution and the computational complexity of the model based on SAA results. It is worth noting that since the computational complexity is increased at least by $|S|$  as the sample size increases, so, we prefer to select a small sample size with a big replication for the model.

We did the sensitivity analysis by considering sample sizes $|S|\in \{ 10,15,20,25,30\} $, number of replication $|M|\in \{ 10,20,30,40,50\} $ and setting the reference sample size to $|S^{'} |=500$. To perform this analysis, the SAA scheme is applied to a 5-node network where two popular random distribution functions namely uniform and normal distribution function are assumed for uncertain network costs in which the non-negativity of the parameter values is preserved by truncating random distributions.

\begin{figure}[h!]
\centering
\includegraphics[scale=0.35]{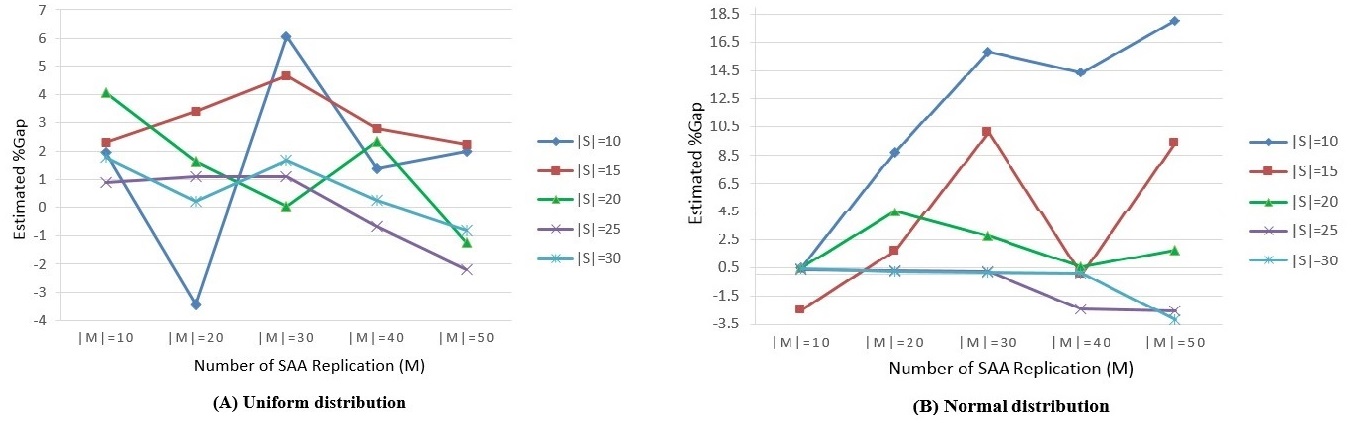}
\caption{\% Gap for a 5-node network}\label{fig2}
\end{figure}

 In the follow, a sensitivity analyses is used to select an optimal sample size. The trends of estimated percent gaps for different sample sizes $|S|$ and replications $|M|$ are show in Fig \ref{fig2}. 

Figure \ref{fig2} demonstrate the trend of estimated percent gap for sample sizes 10, 15, 20 with increasing replication $|M|$ for both random distribution functions, while it uniformly decreased for sample size 25 and never exceeds \%1. On the other hand the behavior of sample size 25 is better when compared with sample size 30 except for $|M|=20$ in Figure \ref{fig2}A and $|M|=50$ in Figure \ref{fig2}B. It is possible that the randomness of parameters led to this deviation. Furthermore, comparing our results for different replications show that by increasing the number of replication, the estimated percent gap is decreased.

\begin{figure}[h!]
\centering
\includegraphics[scale=.34]{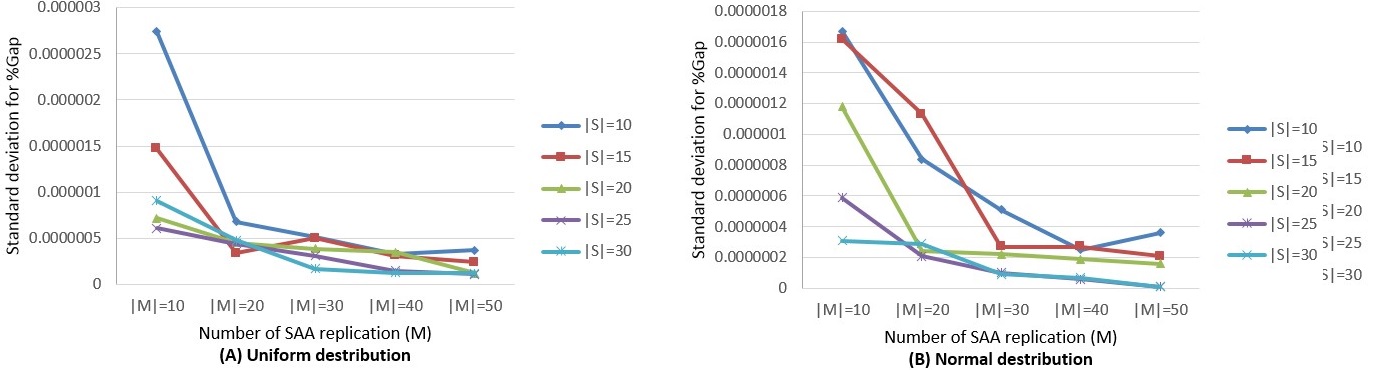}
\caption{Standard deviation for \% Gap for a 5-node network}\label{fig3}
\end{figure}

Figure \ref{fig3} plot the estimated standard deviation for the gap for different sample sizes $|S|$ and replications $|M|$.

Figure \ref{fig3} shows the trends of the standard deviation for the estimated optimality gap for different sample sizes and replications for both random distribution functions. The standard deviation generally is decreased by increasing sample size and number of replications. Also, for sample sizes 25 and 30 the deviation is significantly less than other sample sizes for each replication, especially in Figure \ref{fig3}B. Comparing the standard deviation for sample size 25 with sample size 30, we can say that the sample size 25 is more accurate than sample size 30 except for $|M|=30$, in Figure \ref{fig3}A and $|M|=10$, in Figure \ref{fig3}B.

\begin{figure}[h!]
\centering
\includegraphics[scale=0.35]{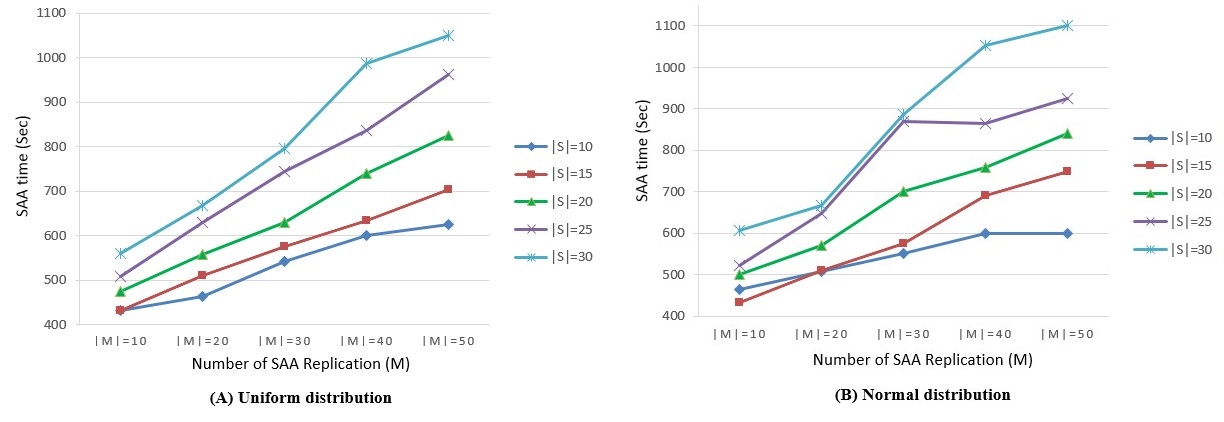}
\caption{Total CPU time for a 5-node network}\label{fig4}
\end{figure}

Figure \ref{fig4}A and B show that the required CPU time for solving the SAA problem for different sample sizes and replications numbers $|M|$ for the uniform and normal distribution functions, respectively.

Figure \ref{fig4} indicate that the required CPU time for solving the SAA problem is nearly linearly increasing with increasing problem size in both cases. For example, the required CPU time for sample size $|S|=10$ and replication $|M|=10$ is 432 seconds whereas it only increases to 625 seconds for replication $|M|=50$. We see that the sample size $|S|=25$ seems be an optimal sample size that have a good trend of the optimal gap as well as the required CPU time. Consequently, we carried out the rest of computational experiments assuming a sample size of 25.

In the end, a comparison measure, namely the Value of the Stochastic Solution (VSS), is given to quantify the necessity of the uncertainty analysis (Birge and Louveaux, 2011). VSS is formulated as the difference between expectation of the expected value problem (EEV) and the average cost values of the problems with random variables (RP). For a specific problem with 25 scenarios, the EEV is 468.2 where the RP is obtained as 490.59. Hence a \%4.6 cost can be saved due to the consideration of uncertainty generated by the random variables. 

\subsection{The effect of risk consideration in the hub location model}

\begin{figure}[h!]
\centering
\includegraphics[scale=0.28]{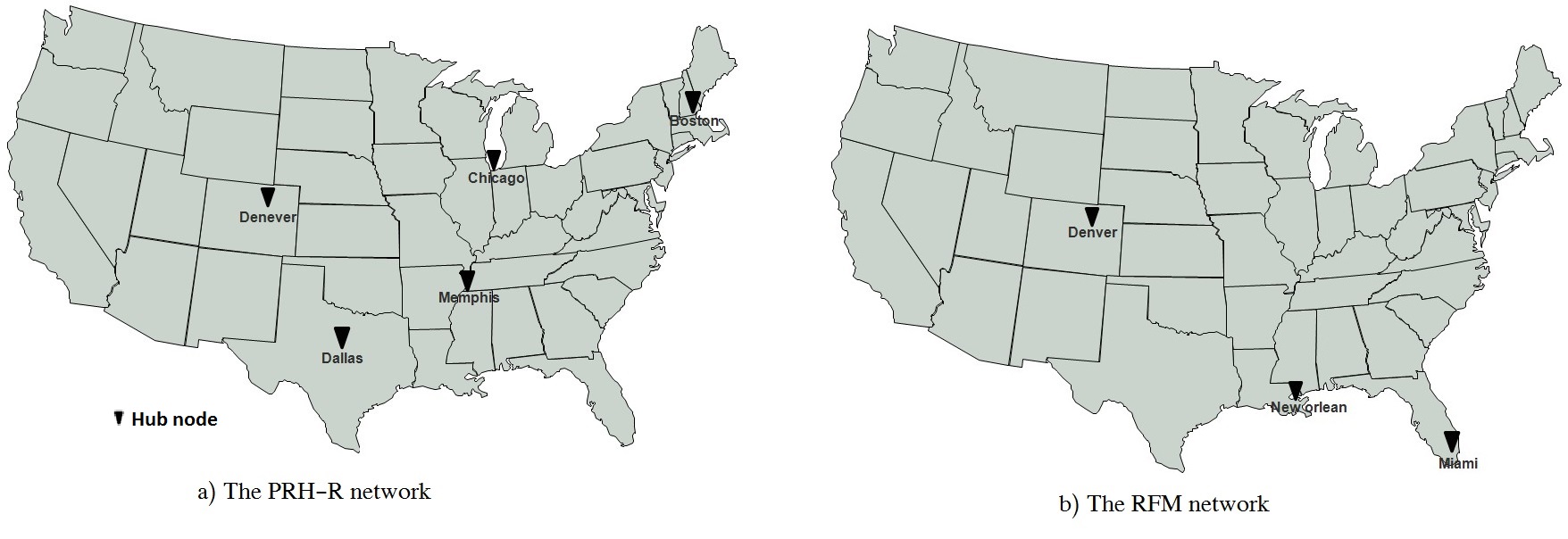}
\caption{Comparison of output for PRH-R and RFM networks}\label{fig5}
\end{figure}

A comparison experiment between the PRH-R and risk free models (RFM) is implemented to show the effectiveness of the proposed approach. RFM, as implied, does not take into account any risk parameters. The hub network related the first period time and scenario is represented for two networks in Figure \ref{fig5}.

The effectiveness of risk consideration can be seen from the obtained results. The set of hubs for PRH-R model are \{Denever, Dallas, Memphis, Chicago, Boston\}, while for the RFM the set of hubs including \{Denever, New orlean, Miami\}. Moreover, the greater number of located facilities in the PRH-R network confirms that the possibility of failure that is taken into consideration in the proposed robust model as compared with the classical version. Furthermore, 1000 failure scenarios were simulated with failure probability for a predetermined $i-j$ flows. The lost flows for both the RFM and PRH-R networks were calculated as a comparison criterion. The obtained results show that there are 9117 unserved commodities for the PRH-R network, while we will would have faced over 12751 unserved commodities in the RFM network (39\% more). 

\subsection{Proposed heuristic performance}

In the following, the performance of the proposed relax-and-decomposition method is evaluated by comparing the results obtained from the solution of the proposed model using a classical LR algorithm and GAMS software. It is also important to show the MPBD justification. So, a set of computational experiments is implemented to confirm the capability of improved BD method for solving large instances.

\subsubsection{Performance evaluation of the relax-and-decomposition method}

Our first computations analyze the efficiency of the proposed relaxed and decomposition algorithm. A set of instances with different sizes is solved by the proposed solution method, classical LR algorithm and an LP relaxation of PRH-R where all integrality requirements are relaxed. The obtained results are presented in Table \ref{ta2}.

\begin{table}[h!]
{\tiny
\caption{Results of different-sized instances for the LP relaxation \& LR \& LRMPBD on PRH-R}
\tabcolsep=3pt
\begin{tabular}{cccccccccccccc}
\hline 
{}
&
{}
& \multicolumn{2}{p{0.9in}}{\centering \# variables} & \multicolumn{2}{p{0.7in}}{\centering \# constraints} & \multicolumn{2}{p{0.7in}}{\centering LP solution} & \multicolumn{3}{p{1.3in}}{\centering LR solution} & \multicolumn{3}{p{1.7in}}{\centering LRMPBD solution}   \\ \hline 
\textbar H\textbar
& 
\textbar T\textbar
&
Continues & Binary &Equality & Inequality & ${lB}_{LP}$\textit{}  & CPU${}_{(s)}$ & ${lB}_{LR}$\textit{} & CPU${}_{(s)}$\textit{} & \%Dev\textit{${}_{Dif-LR}$} & ${lB}_{{\rm LR-MPBD}}$\textit{} &  CPU${}_{(s)}$\textit{} &  \%Dev\textit{${}_{Dif-LR}$}
\\ \hline
\cr

6
&
6
&
5401
&
32616
&
150
&
22597
&
180.3
&
11
&
59.9
&
15
&
65.6
&
174.6
&
993
&
0
\cr
10
&
10
&
25001
&
251000
&
250
&
102725
&
527
&
286
&
129.2
&
1061
&
73.8
&
493.4
&
1740
&
0
\cr
12
&
12
&
43201
&
520128
&
300
&
176713
&
754.4
&
334
&
551.8
&
902
&
23.2
&
717.3
&
2074
&
0
\cr
14
&
14
&
68601
&
963144
&
350
&
279717
&
1520.6
&
899
&
1520.9
&
936
&
0
&
1460.7
&
3679
&
3
\cr
16
&
16
&
102401
&
1642496
&
400
&
416537
&
135.7
&
960
&
126.7
&
1158
&
0
&
113.6
&
3560
&
10
\cr
18
&
18
&
145801
&
2630232
&
450
&
591973
&
208.4
&
1435
&
185.0
&
1225
&
0
&
147.3
&
4638
&
20
\cr
20
&
20
&
200001
&
4008000
&
500
&
810825
&
1.4
&
1080
&
2451.5-
&
1272
&
99.4
&
14.7-
&
2441
&
0
\cr
20
&
80
&
800001
&
64128000
&
2000
&
3243225
&
9248.1
&
1385
&
-
&
1-
&
-
&
-1.3
&
1490
&
-
\cr
20
&
100
&
1000001
&
1E+08
&
2500
&
4054025
&
10221.6
&
1491
&
-
&
1-
&
-
&
1149
&
7112
&
-
\cr
20
&
120
&
1200001
&
1/44E+08
&
3000
&
4864825
&
15513.6
&
3021
&
-
&
1-
&
-
&
15451
&
8195
&
-
\cr
20
&
130
&
1300001
&
1/69E+08
&
3250
&
5270225
&
-
&
2-
&
-
&
1-
&
-
&
13801.6
&
8939
&
-
\cr
20
&
140
&
1400001
&
1/96E+08
&
3500
&
5675625
&
-
&
2-
&
-
&
1-
&
-
&
15916.7
&
7423
&
-
\cr
20
&
150
&
1500001
&
2/25E+08
&
3750
&
6081025
&
-
&
2-
&
-
&
1-
&
-
&
-11.4
&
8438
&
-
\cr
20
&
160
&
1600001
&
2/57E+08
&
4000
&
6486425
&
-
&
2-
&
-
&
1-
&
-
&
-17.3
&
8266
&
-
\cr
20
&
170
&
1700001
&
2/9E+08
&
4250
&
6891825
&
-
&
2-
&
-
&
2-
&
-
&
-16.7
&
6530
&
-
\cr
20
&
180
&
1800001
&
3/25E+08
&
4500
&
7297225
&
-
&
2-
&
-
&
2-
&
-
&
-28.8
&
7382
&
-
\cr
20
&
190
&
1900001
&
3/62E+08
&
4750
&
7702625
&
-
&
2-
&
-
&
2-
&
-
&
-25.3
&
8963
&
-
\cr
20
&
200
&
2000001
&
4/01E+08
&
5000
&
8108025
&
-
&
2-
&
-
&
2-
&
-
&
-1.34
&
9865
&
-
\cr

\multicolumn{14}{|p{1in}|}{${}^{1}$- : Time limitation\textit{\newline }${}^{2}$- : Lack of memory} \\ \hline 
\end{tabular}\label{ta2}}
\end{table}

According to the Table \ref{ta2}, efficiency of the proposed relax-and-decomposition algorithm is remarkable, especially for the large scale instances. In small size instances, the proposed solution method can find a lower bound for the proposed model as well as the classical LR algorithm. But, we can see that for sizes $T=\{6, 10, 12, 20\}$, a more efficient lower bound can be obtained by the our relaxed and decomposition method. Recall that a lower bound is achieved for instances greater than T=120 by our solution method while time and memory limitation led to no solution for the classical LR and LP relaxation problems. In particular, LRMPBD algorithm could find a more efficient lower bound in comparison with the classical LR algorithm in 57\% of these instances and for two the algorithm could reach the solution. Furthermore, the percent gap of LRMPBD algorithm never exceeds 20\%. These experiments confirm that the proposed solution method outperforms the classical LR algorithm and LP relaxation problem.

\subsubsection{Performance evaluation of the improved Benders decomposition algorithm}

We presented a set of numerical results to investigate the behaviour of MPBD method. The first round of experiments evaluated the capability of multi-Pareto cuts comparing with other different cut generators as single-cut, multi-cut, Pareto-cut and multi-Pareto cut versions to find the solution. The obtained results are presented in Table \ref{ta3}.
\begin{table}[h!]
\caption{Results of different-sized instances for different BD versions $\left(\times {10}^2\right)$}\label{ta3}
\centering
{\tiny
\begin{tabular}{|p{0.12in}|p{0.34in}|p{0.33in}|p{0.35in}|p{0.36in}|p{0.23in}|p{0.3in}|p{0.3in}|p{0.3in}|p{0.3in}|p{0.3in}|p{0.3in}|p{0.3in}|p{0.1in}|p{0.3in}|} \hline 
\multicolumn{1}{p{0.12in}}{}&\multicolumn{2}{p{0.7in}}{Number of Variables} &\multicolumn{2}{p{0.7in}}{Number of constraints} &\multicolumn{2}{p{0.7in}}{SBD solution} &\multicolumn{2}{p{0.7in}}{MBD solution} &\multicolumn{2}{p{0.7in}}{PBD solution} &\multicolumn{2}{p{0.7in}}{MPBD solution} \cr  \hline
 \textbar T\textbar
 &
 Continues
 &
 Binary
 &
 Equality
 &
 Inequality
 &
 CPU(s)
 &
 Dev(BD)
 &
 CPU(s)
 &
 Dev(BD) 
 &
 CPU(s)
 &
 Dev(BD) 
 &
 CPU(s)
 &
 Dev(BD)
 \cr \hline
10
&
100001
&
1002000
&
250
&
405425
&
1178
&
1.69665
&
1304
&
1.50967
&
1192
&
0.87619
&
1523
&
0
\cr
20
&
200001
&
4008000
&
500
&
810825
&
1354
&
0.17258
&
2050
&
0.13981
&
1799
&
0.00246
&
2163
&
0
\cr
30
&
300001
&
9018000
&
750
&
1216225
&
1909
&
0.47578
&
2504
&
0.32280
&
2046
&
0
&
2689
&
0.13920
\cr
40
&
400001
&
16032000
&
1000
&
1621625
&
2965
&
4.81534
&
3536
&
4.76652
&
2813
&
0.56201
&
3678
&
0
\cr
50
&
500001
&
25050000
&
1250
&
2027025
&
3624
&
1.56830
&
5532
&
1.559976
&
4800
&
0
&
5162
&
0.04208
\cr
60
&
600001
&
36072000
&
1500
&
2432425
&
5020
&
19.9793
&
6463
&
19.28069
&
5958
&
1.56265
&
6445
&
0
\cr
70
&
700001
&
49098000
&
1750
&
2837825
&
6443
&
0.47404
&
8137
&
0.468368
&
7282
&
0.11495
&
6637
&
0
\cr
80
&
800001
&
64128000
&
2000
&
3243225
&
7094
&
0.04190
&
8751
&
0.039082
&
8151
&
0
&
7732
&
0.03796
\cr
90
&
900001
&
81162000
&
2250
&
3648625
&
8798
&
0.91757
&
9421
&
0.002854
&
9293
&
0.91590
&
9372
&
0
\cr
100
&
1000001
&
1E+08
&
2500
&
4054025
&
8664
&
0.09882
&
9692
&
0
&
9009
&
0.09883
&
9445
&
0
\cr
120
&
1200001
&
1/44E+08
&
3000
&
4864825
&
1-
&
-
&
1-
&
-
&
1-
&
-
&
9640
&
-
\cr
0
&
130
&
1300001
&
169338000      3250
&
5270225
&
1-
&
-
&
2-
&
-
&
2-
&
-
&
9768
&
-
\cr
0
&
140
&
1400001
&
196392000      3500
&
5675625
&
2-
&
-
&
2-
&
-
&
2-
&
-
&
9978
&
-
\cr
0
&
150
&
1500001
&
225450000      3750
&
6081025
&
2-
&
-
&
2-
&
-
&
2-
&
-
&
9856
&
-
\cr
Avg.
&
-
&
-
&
-
&
-
&
-
&
3.02404
&
-
&
2.81165
&
-
&
0.41330
&
-
&
0.02192
\cr
\hline

\end{tabular}
}
\end{table}

The effect of employing the multi-Pareto strategy to improve the quality of the solution can be seen from the results presented in Table 3. According to these results it can be concluded that the relative average gap of multi-Pareto strategy is about \%99 less than the value corresponding to other strategies. Moreover, on average there is just a \%16 increase in CPU time due to the more effort expended in generated cuts in each iteration compared with other improvement strategies, where the CPU time increase is always bellow \%37. For larger instances in which multi-cut and Pareto-cut strategies are unable to reach a solution in a predefined CPU time, the multi-Pareto strategy is shown to converged in a reasonable time. This confirms the superiority of our solution algorithm in solving the PRH-R problem.

\begin{figure}[h!]
\centering
\includegraphics[scale=0.35]{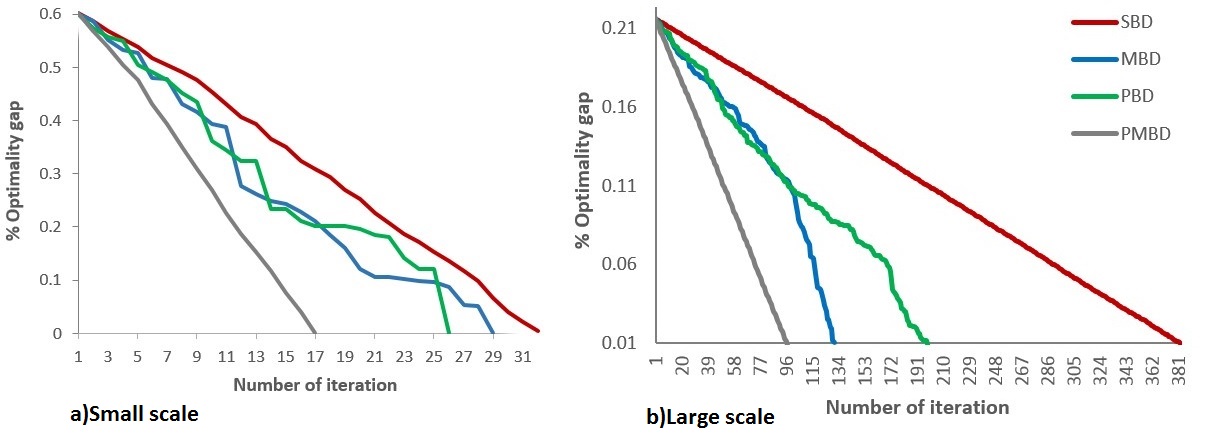}
\caption{Behavior of optimality gap for different cutting strategies during different iterations}\label{fig6}
\end{figure}

 The trend for convergence for different cut generator procedures are shown in Figure \ref{fig6}. The performance improvement is demonstrated for small and large size instances.

The convergence of multi-Pareto strategy as well as other cut generators is investigated in Figure \ref{fig6}. For small size, the multi-Pareto procedure converged in 17 iterations while the single cut strategy requires 33 iterations. Similarly, the multi-Pareto strategy caused a 286 decreased in iterations as compared to a single cut generator for large scale instances. We can see that increasing the size of problem led to more complexity in the structure of RMP constraints and thus resulted in the problem being time consuming. In this situation, adding multi-cuts in each iteration can restrict the solution space of the RMP, much more effectively than other approaches.

\paragraph{Core point update }

\begin{figure}[h!]
\centering
\includegraphics[scale=0.35]{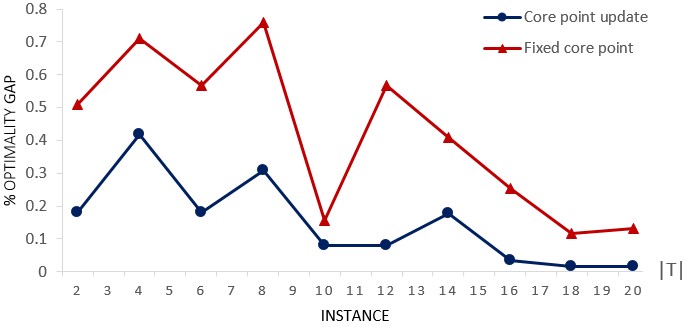}
\caption{The impact of the core point updating rule}\label{fig7}
\end{figure}

We use a linearization combination procedure as the update rule for the core points. Now, a comparison experiment is implemented to show the effectiveness of this rule in finding the optimal solution. Results in Figure \ref{fig7} confirm that BD algorithm based on this updating rule is much more efficient.

\section{Conclusion}

In this research, a nonlinear $p$-robust hub location problem was proposed. An augmented chance constraint was applied in a min-max regret form. Also, some unpredictable factors such as congestion, security, delay time and regional air pollution were considered as risk factors. The problem was formulated as a scenario based model. In this way, a sampling method, namely a sample average approximation was applied. Then, a relax-and-decompose heuristic was shown to compute a good upper bound for the proposed model where a multi-Pareto cut BD method was incorporated into a LR algorithm in order to provide strong lower bounds for the large-scale instances. The capability of the proposed model as a robust network design was verified via numerical experiments. Finally, a sensitivity analysis demonstrated the accuracy of the solution procedure for solving the proposed model as compared with the classical Lagrangian relaxation method. For future research, employing branching methods in modeling the master problem may be more effective as an alternative accelerating strategy to solve the PRH-R problem. Further research will include the applications of the proposed solution method to other well-known hard problems such as hierarchical HLPs or vehicle routing problems.

\section*{References}

\bibliography{mybibfile}

\noindent URL: \url{https://en.m.wikipedia.org/wiki/Valigonda_train_wreck}, last visited: June 2016. 

\noindent URL: \url{https://en.wikipedia.org/wiki/Air\_travel_disruption_after_the_2010_Eyjafjallaj%C3%B6kull_eruption}, last visited: June 2016. 

\noindent URL:\url{ https://en.wikipedia.org/wiki/Effects_of_Hurricane_Katrina_in_New_Orleans}, last visited June 2016. 

\citep{1,2,3,4,5,6,7,8,9,10,11,12,13,14,15,16,17,18,19,20,21,22,23,24,25,26,27,28,29,30,31,32,33,34,35,36,37,38,39,40,41,42,43,44,45,46,47,48,49,50,51,52,53,54,55,56,57,58,59,60,61,62,63,64,65,66,67,68,69}
\end{document}